\documentclass[12pt]{amsart}
\usepackage{amsfonts, amsmath}
\usepackage{ amsmath, amsthm, amsfonts, amssymb, color, graphics, graphicx}
 \usepackage{mathrsfs}
 \usepackage{amsmath,amstext,amsthm,amssymb,amsxtra,verbatim}
\newtheorem{theorem}{Theorem}[section]
\newtheorem{proposition}[theorem]{Proposition}

\newtheorem{lemma}[theorem]{Lemma}
\newtheorem{definition}[theorem]{Definition}

\newtheorem{remark}[theorem]{Remark}

\newcommand\D{\mathcal{D}}

\newcommand\R{\mathbb{R}}

\newcommand\ZZ{\mathbb{Z}}

\def \l {\lambda}
\newcommand{\eps}{\varepsilon}

\newcommand{\supp}{{\rm supp}{\hspace{.05cm}}}

\title[The bilinear Bochner-Riesz  problem]
{The bilinear Bochner-Riesz  problem}

\medskip
\author [  F. Bernicot,    L. Grafakos,    L. Song    and    L.X. Yan    ]
 {Fr\'ed\'eric Bernicot,     Loukas Grafakos,    Liang Song     and     Lixin Yan    }

\subjclass[2010]{42B15, 42B20,   47B38.}
\keywords{Bilinear Bochner-Riesz  problem;  bilinear Fourier multipliers;  bilinear restriction-extension operators}

\date{\today}

\begin{document}

\begin{abstract}
Motivated by   the problem of spherical summability of products of Fourier series, we study the boundedness
of the bilinear Bochner-Riesz multipliers $(1-|\xi|^2-|\eta|^2)^\delta_+$ and we make some advances in this
investigation.  We obtain an optimal result concerning the boundedness of these means from
$L^2\times L^2 $ into $L^1$ with minimal smoothness, i.e., any $\delta>0$,
   and we   obtain estimates for other pairs of spaces
for larger values of $\delta$. Our study is broad enough   to encompass   general
 bilinear multipliers $m(\xi,\eta)$ radial in $\xi$ and $\eta$
 with minimal smoothness,    measured in Sobolev space norms.
The  results obtained 
are based on a variety of techniques, that include  Fourier series expansions, orthogonality,  and
bilinear restriction and extension theorems.
\end{abstract}

\maketitle

 \tableofcontents

%\bigskip

\section{Introduction  }
\setcounter{equation}{0}

\medskip

The study of the summability of the product  of two $n$-dimensional Fourier series  leads to questions
concerning  the  norm convergence  of  partial sums  of the form
$$
\sum_{|m|^2+|k|^2\le R^2} \widehat{F}(m) e^{2\pi i m\cdot x} \,\widehat{G}(k) e^{2\pi i k\cdot x},
$$
as $R\to \infty$, or more generally, of the bilinear Bochner-Riesz means
\begin{equation}\label{bBBR}
\sum_{|m|^2+|k|^2\le R^2} \big(1-\tfrac{|m|^2+|k|^2}{R^2}\big)^\delta \widehat{F}(m) e^{2\pi i m\cdot x}\,
 \widehat{G}(k) e^{2\pi i k\cdot x}
\end{equation}
for some $\delta\ge 0$. Here $F,G$ are $1$-periodic functions on the $n$-torus and $\widehat{F}(m), \widehat{G}(k)$
are their Fourier coefficients and $m,k\in \mathbb Z^n$.
The bilinear   Bochner-Riesz problem is the study of the norm convergence of the sum in \eqref{bBBR}.
By basic functional analysis and transference, this problem is equivalent to the study of the
$L^{p_1}\times L^{p_2}\to L^p$ boundedness of the
bilinear Fourier multiplier operator
\begin{equation}\label{bBBR2}
 S^{\delta} (f, g)(x)  := \iint_{|\xi|^2+|\eta|^2\leq 1}
  \big(1-  |\xi|^2-|\eta|^2\big)^{\delta} {\widehat f}(\xi) {\widehat g}(\eta)e^{2\pi ix\cdot (\xi+\eta)} d\xi d\eta  .
  \end{equation}
Here $x\in \R^n$, $f,g$ are functions   on $\R^n$ and ${\widehat f},  {\widehat g}$ are their Fourier transforms.

The Bochner-Riesz summability question is a fundamental problem in mathematics. Its study has led to the
development of important notions, tools, and results in Fourier  analysis,  and has created numerous directions of
 research. The Bochner-Riesz conjecture is well known to be difficult and remains unsolved for   indices $p$
  near $2$ in dimensions $n\ge 3$. The  bilinear   Bochner-Riesz problem
is   more difficult than its linear counterpart because of the natural complexity that arises from the
mixed summability and also from the shortage of
 techniques to study   bilinear Fourier multipliers with minimal smoothness.
 The present work is motivated by this problem and
 fits   under the scope of the   program to find minimal smoothness conditions for a bilinear
Fourier multiplier to be bounded on products of Lebesgue spaces.
We are mainly interested in theorems concerning  compactly supported Fourier multipliers.
 The main question we   address  is what is the least amount of differentiability
 required of a generic function   on $\R^n\times \R^n$ to become a bilinear Fourier multiplier on a certain
 product of Lebesgue spaces. For the purposes of this article,  differentiability is measured in terms of
Sobolev space norms which   quantitatively fine-tune
 fractional smoothness. Our results concerning the bilinear Bochner-Riesz means
 fit in this  general
 framework.

 It is well known that linear multiplier operators are $L^2$ bounded if and only if the multiplier is a bounded function.
But we know  from \cite{GK2} that there exist  smooth functions $m$ satisfying
$$
|\partial^{\alpha}_{\xi} \partial^{\beta}_{\eta} m(\xi, \eta)|\leq C_{\alpha \beta} |\xi|^{-|\alpha|} |\eta|^{-|\beta|}\ , \qquad
\xi,\eta\neq 0
 $$
for all multi-indices $\alpha,\beta$ and also from \cite{BT} that there exist  smooth functions $m$ satisfying
$$
|\partial^{\alpha}_{\xi} \partial^{\beta}_{\eta} m(\xi, \eta)|\leq C_{\alpha \beta}
 $$
for all $(\xi, \eta)\in \R^{2n}$ and all multi-indices $\alpha$ and $\beta$,
 which do not give rise to bounded bilinear operators (as defined in \eqref{multdef})
 from   $L^{p_1}(\R^n)\times L^{p_2}(\R^n)$ to $L^p(\R^n)$ when
$1/p_1 + 1/p_2=1/p$  and $1\leq p_1, p_2, p\leq \infty.$ So there is no direct analogy with the linear case where
$L^2$  presents itself as a natural starting point of the investigation of multiplier theorems.

So we aim to focus our study on more particular bilinear operators. Suppose that a bilinear operator $T$,
initially acting from $\mathscr  S(\R^n) \times \mathscr  S(\R^n)$ to $\mathscr S'(\R^n)$, admits an $L^{p_1}\times L^{p_2}\to L^p$
bounded extension, i.e., it is a bilinear Fourier multiplier
for some $1<p_1,p_2,p<\infty$ with $1/p_1+1/p_2=1/p$. Then the following properties are equivalent:

\smallskip
\noindent (i) {\it Frequency representation}. There exists a bounded function $m$ on $\R^{2n}$ such that for all $f,g,h\in \mathscr  S(\R^n)$ we have
$$
\int_{\R^n} T(f,g)(x) \overline{h(x)}\, dx = \int_{\R^{2n}} m(\xi,\eta) \widehat{f}(\xi) \widehat{g}(\eta)
\overline{\widehat{h}(\xi+\eta)}\,  d\xi d\eta.
$$

\noindent (ii) {\it Kernel representation}. There exists a tempered distribution $K$ on $\R^{2n}$
such that for all $f,g  \in \mathscr  S(\R^n)$ we have
$$
 T(f,g)(x) = \langle K , f(x-\cdot) \otimes g(x-\cdot) \rangle ,
$$
where $(f(x-\cdot) \otimes g(x-\cdot) )(y,z)= f(x-y)g(x-z)$ for all $x,y,z\in \R^n$.

\noindent (iii)  {\it Commutativity with simultaneous translation}.  For every $y\in \R^n$ and for every function $f,g\in
 \mathscr  S(\R^n)$ we have
$$ T(\tau_y(f), \tau_y(g)) = \tau_y(T(f,g))$$
where $\tau_y$ is the translation operator $\tau_y(f)(x)=f(x-y)$.
 This property takes into account the additive structure of the Euclidean space via the group
of translations.

Bilinear multipliers are not invariant under rotations but the following is true:
let $T$ be a bilinear Fourier multiplier on $\R^n$ and $m$ be its symbol; then the symbol is biradial, i.e.,
 $m(\xi,\eta) = m_0(|\xi|,|\eta|)$ (for some $m_0\in L^\infty(\R^2)$) if and only if for every  pair of
orthogonal transformations (rotations)  ${\mathcal R}_1,\ {\mathcal R}_2 $ of $\R^n$ we have
$$ T(f,g)(0) = T(f\circ {\mathcal R}_1, g\circ {\mathcal R}_2)(0).$$
Such operators naturally appear  in the study of scattering properties associated to quadratic PDEs
involving functions of the Laplacian   (see \cite[Section 2.3]{BBMNT}).

Of course, this property reduces, in some sense, a $2n$-dimensional symbol to a $2$-dimensional symbol
and this work aims to understand how one can take advantage of this property.
We observe that for  radial  multipliers,  differentiability is only relevant  in the  radial direction,
and the point $L^2\times L^2 \to L^1$ seems to be the one requiring
 the least smoothness. We point out that the duals of a bi-radial bilinear multiplier $m_0(|\xi|,|\eta|)$,
 $m_0(|\xi+\eta|,|\eta|)$  and  $m_0(|\xi|, |\xi+\eta |)$,    are not bi-radial functions, so certain results
 we obtain are not symmetric in the local $L^2$ triangle, i.e., the   set $\{(1/p_1,1/p_2,1/p)$ with $2\le p_1,p_2,p'\le \infty\}$;   here
 $p'=p/(p-1)$.

Let us give some examples of bilinear multipliers, pointing out different situations with respect to the nature
of the singular space of the symbol $m$: we say that $m$ is allowed to be singular a set $\Gamma \subset \R^{2n}$
if $m$ is smooth in the complement $\Gamma^c$ and satisfies
\begin{equation} \left| \partial_{(\xi,\eta)}^\alpha m(\xi,\eta) \right|
\leq C_{\alpha} d((\xi,\eta),\Gamma)^{-|\alpha|} \label{eq:sing}
\end{equation}
for every $(\xi,\eta)\in \Gamma^c$ and multi-index $\alpha$.

\smallskip
\begin{itemize}
\item {\it Singularity at one point $\Gamma:=\{0\}$} (Coifman and Meyer \cite{CM1, CM2, CM3}.)
Suppose that the bounded function $m(\xi,\eta)$ on $\R^{2n}$ satisfies (\ref{eq:sing}) with $\Gamma:=\{0\}$
and so $d((\xi,\eta),\Gamma) \simeq |\xi|+|\eta|$. Then the operator $T_m$ is bounded from
$L^{p_1}(\R)\times L^{p_2}(\R)$ to $L^p(\R)$ when $1/p_1+1/p_2=1/p$, $1<p_1,p_2,p\le\infty$.
This theorem was extended to the case $1/2<p\le 1$ by Grafakos and Torres \cite{GT} and independently
 by Kenig and Stein \cite{KS}. This extension also includes the endpoint case $L^1\times L^1 \to
L^{1/2,\infty}$.

\item {\it Singularity along a line} (Lacey and Thiele \cite{LT1,LT2}.) The bilinear Hilbert
transform was shown to be bounded on Lebesgue spaces by Lacey and Thiele. This corresponds to the case where $\Gamma$
is a non-degenerate line of $\R^2$.

\item {\it Singularity along the circle, $\Gamma:={\mathbb S}^1$} (Grafakos and Li \cite{GL}.)
The characteristic function of the unit disc
is a bilinear Fourier multiplier from $L^{p_1}(\R)\times L^{p_2}(\R)$ to  $L^p(\R )$ when  $2\leq p_1, p_2,p'<\infty$ and $1/p_1 + 1/p_2=1/p$.

\item {\it Singularity on the boundary of a disc, $\Gamma:={\mathbb S}^1$} (Diestel and Grafakos \cite{DG}.)
 The characteristic function of the unit disc
in $\R^4$ is not a bilinear Fourier multiplier  from $L^{p_1}(\R^2)\times L^{p_2}(\R^2)$  to $L^p(\R^2)$  when
$1/p_1 + 1/p_2=1/p$  and exactly one of $p_1, p_2, p'$ is less than $2$.

\item {\it Singularity along a curve} (Bernicot-Germain \cite{BG2}.) In this work, certain one-dimensional
 bilinear operators  whose symbols are singular along a curve are shown to be bounded. Taking advantage of
 the non-degeneracy  or the non-vanishing curvature some sharp estimates in the H\"older scaling
 (or sub-H\"older scaling) are proved. There, the
variables are uni-dimensional and $\Gamma$ is a curve in $\R^2$ and so it has dimension $1$.

\item {\it Singularity along a subspace} (Demeter, Pramanik and Thiele \cite{DT,DPT}.) In \cite{DPT}, if $\Gamma$ is a subspace,
preserving the ``$n$-coordinates structure'' and of dimension $\kappa \leq \frac{3d}{2}$, then operators associated to symbols
singular along such non-degenerate subspace are shown to be bounded on Lebesgue spaces \cite{DPT}. However, the time-frequency
analysis used for the bilinear Hilbert transform is not adapted to the multi-dimensional setting with a high-dimensional singular
 subspace (as observed in \cite{DT}). Indeed, it does not allow to understand how the mixing of the coordinates behave in the
 frequency plane. A simpler model was considered by Bernicot and Kovac to handle   the ``twisted paraproducts" \cite{B,K}.

\item {\it Boundedness on Hardy spaces} (Miyachi and Tomita \cite{MT}, Tomita  \cite{TN}.) Suppose that
 $0<p_1, p_2\leq \infty$ and ${1\over p_1}+{1\over p_2}={1\over p}$ and that
$$
s_1>\max \Big\{ {n\over 2}, \ {n\over p_1}-{n\over 2}\Big\}, \ \ \  s_2>\max \Big\{ {n\over 2}, \ {n\over p_2}-{n\over 2} \Big\},$$
 $$
s_1 +s_2 >n\Big( {1\over p_1} + {1\over p_2} -{1\over 2}\Big) =n\Big({1\over p} -{1\over 2}\Big).$$
Assume that for some smooth bump $\Psi$ supported in $6/7\le |\xi|\le 2$ and equal to $1$ on $1\le |\xi |\le 12/7$
we have
$$
K=\sup_{j\in \mathbb Z} \| m(2^j \xi_1, 2^j \xi_2) \Psi (\xi_1, \xi_2) \|_{W^{(s_1, s_2)}}<\infty,
$$
where
$$
\| F\|_{W^{(s_1,s_2)}} = \Big( \int_{\R^{2n}} (1+|\xi_1|^2)^{2s_1} (1+|\xi_2|^2)^{2s_2} |\widehat F(\xi_1,\xi_2)|^2
 d \xi_1d\xi_2 \Big)^{1/2}.
$$
Then $T_m$ is a bounded bilinear operator on products of Hardy spaces with norm
$$
\|T_m\|_{H^{p_1}(\R^n)\times H^{p_2}(\R^n) \to L^p({\R^n})}\leq C \, K  ,
$$
where $L^\infty({\R^n})$ should be replaced by $ BMO({\R^n}) $ when $p_1=p_2=\infty$.

 \end{itemize}

From this quick review of   existing results, it appears that high-dimensional symbols singular along
hypersurfaces have not been studied, according to our understanding.
Our approach of biradial bilinear Fourier multiplier will allow us to consider the bilinear counterpart
$S^\delta$ (as defined in \eqref{bBBR2})
of the celebrated Bochner-Riesz multiplier.
%$$ S^{\delta} (f, g)(x)  := \iint_{|\xi|^2+|\eta|^2\leq 1} e^{2\pi ix\cdot (\xi+\eta)}
%  \big(1-  |\xi|^2-|\eta|^2\big)^{\delta} {\widehat f}(\xi) {\widehat g}(\eta)d\xi d\eta. $$
Here the symbol is singular along the sphere $\{ (\xi,\eta)\in\R^{2n},\ |\xi|^2+ |\eta|^2 =1\}$ which has dimension $2n-1$.

The almost optimal solution of the bilinear Bochner-Riesz problem in dimension $1$ is outlined in Theorem
\ref{BR-1dim}.
We end this introduction by summarizing some critical estimates obtained in this article for the bilinear
Bochner-Riesz means when $n\ge 2$:

\begin{itemize}

\item       $S^{\delta}$ is bounded from
 $L^2(\R^n)\times L^2(\R^n)$ to $L^1(\R^n)$ when    $\delta>0$.  (Theorem \ref{thm:221}.)

 \item   $S^{\delta}$ is bounded from
 $L^2(\R^n)\times L^\infty(\R^n)$  to $L^2(\R^n)$ if $\delta>\frac{n-1}{2}$. (Theorem \ref{th4.7}.)

\item   $S^{\delta}$ is bounded from
 $L^{1}(\R^n)\times L^{\infty}(\R^n)$ to $L^1(\R^n)$ when  $\delta>{n\over 2}$. (Theorem \ref{thm4.88}.)

\item $S^{\delta}$ is  bounded from $L^{p_1}({\R^n})\times
L^{p_2}(\R^n)$  to $L^p(\R^n)$ when
$1\leq p_1, p_2 < {2n/(n+1)}$,
 $ {1/p}={1/p_1}+ {1/p_2}$,
$\delta>n\alpha(p_1, p_2)-1,$  where   $\alpha(p_1,p_2)$ is as in \eqref{e3.10}. (Theorem \ref{th4.2}.)

\end{itemize}

\medskip

\section{Notation and preliminary results}
\setcounter{equation}{0}

\subsection{Notation}
We introduce the notation  that will be relevant for this paper.
We use $A\lesssim  B$ to denote the statement that $A\leq CB$ for some implicit, universal constant $C$,
and the value of $C$ may change from line to line.
 We denote by $x\cdot y= \sum_j x_jy_j$
the usual dot product of points $x=(x_1,\dots , x_n)$ and $y=(y_1,\dots , y_n)$ in
$\R^n$. We denote by
 ${\mathscr S}(\R^n)$     the Schwartz space  of all rapidly decreasing smooth functions on $\R^n$.
 For a function $f$ in ${\mathscr S}(\R^n)$, we
define the Fourier transform ${\mathscr F}f$ and its inverse
Fourier transform   ${\mathscr F}^{-1}f$   by the formulae
$$
{\mathscr F}f(\xi) ={\widehat f}(\xi)= \int_{\R^n} e^{-2\pi ix\cdot \xi} f(x)\, dx
$$
and
$$
{\mathscr F}^{-1}f(\xi) = {\check f}(\xi)  =  \int_{\R^n} e^{2\pi ix\cdot \xi} f(x)\, dx.
$$

For $1\leq p\leq \infty$, we denote by $p'$ its conjugate exponent, i.e., the unique
number in $[1,\infty]$ such that
  $1/p+1/p'=1$. For $1\le p\le+\infty$, we denote the
norm of a function $f\in L^p({\mathbb R}^n)$ by $\|f\|_p$.
For $s\ge 0$ and $1<p<\infty$, the Sobolev space $W^{s,p}(\R^n)$ is defined as the space of all functions such that
$$
(I-\Delta)^{2/s}(f) = {\mathscr F}^{-1} \Big( (1+4\pi^2 |\xi|^2)^{s/2} {\mathscr F}f(\xi)\Big)
$$
lies in $L^p({\mathbb R}^n)$. In this case we set $\|f\|_{W^{s,p}} = \|(I-\Delta)^{s/2}(f)\|_{L^p}$.

The scalar product in $L^2({\mathbb R}^n)$ is denote by
$$
\langle f,g\rangle = \int_{\R^n} f(x) \overline{g(x)}\, dx\, .
$$
Let $X,Y,Z$ be quasi-normed spaces. If $T$ is a bounded bilinear operator from $
X\times Y$ to $Z$,
we write $\|T\|_{X\times Y\to Z} $ for
the  operator norm of $T$.
Given a  subset $E\subseteq \R^n$, we  denote by  $\chi_E$   the characteristic
function of   $E$ and we denote by
$$
P_Ef(x)=\chi_E(x) f(x)
$$
the ``projection" operator on $E$.

Given a bounded function $m(\xi,\eta)$ on $\R^n\times \R^n$, we denote by   $T_m$ the
 bilinear Fourier multiplier with symbol $m$. This operator is written in the form
\begin{equation}\label{multdef}
T_m(f,g)(x)=  \int_{\R^n}\int_{\R^n} e^{2\pi ix\cdot(\xi+\eta)} m(\xi, \eta) {\widehat f}(\xi) {\widehat g}(\eta)d\xi d\eta
\end{equation}
for Schwartz functions $f,g$. Equivalently, in physical space is given as
$$
T_m(f,g)(x)= {\mathscr F}^{-1} \bigg[\int_{\R^n}   m(\xi-\eta, \eta) {\widehat f}(\xi-\eta) {\widehat g}(\eta)d\eta \bigg](x)
$$
and also as
\begin{eqnarray}\label{uuu}
T_m(f,g)(x)= \int_{\R^n} \int_{\R^n}  {\widehat m}(y-x, z-x) f(y) g(z) dydz\, .
\end{eqnarray}
This is a bilinear translation invariant
operator with kernel $K(y,z)={\widehat m}(-y , -z)$, i.e., it has the form
\begin{eqnarray}\label{uuu66}
T_m(f,g)(x)= \int_{\R^n} \int_{\R^n}  K(x-y,x-z) \, f(y) g(z)\, dydz\, .
\end{eqnarray}

\subsection{Criteria for boundedness of bilinear multipliers}
We begin with the following trivial situation.

 \begin{lemma}\label{le2.1}  Let $1\leq p_1, p_2, p\leq \infty$ and ${1\over p}={1\over p_1}+{1\over p_2}$.
  If the symbol $m(\xi, \eta)$ satisfies
$$
A_1= \int_{\R^n}\int_{\R^n} |{\widehat m}(x,y)|dxdy <\infty,
$$
then $T_m$ maps $L^{p_1}(\R^n) \times L^{p_2}(\R^n) \to L^{p}(\R^n)$ with
$$
\|T_m \|_{L^{p_1}\times L^{p_2}\to L^p}
\leq  A_{1}.
$$
\end{lemma}

 \medskip
 The proof of Lemma~\ref{le2.1} is omitted since it is an easy consequence of Minkowski's integral inequality and H\"older's inequality.

 %\noindent
%{\bf Proof.} Using   \eqref{uuu}, we write
%\begin{eqnarray*}
%  \|T_m(f,g)\|_p&=& \bigg( \int_{\R^n} \bigg| \iint_{\R^{2n}} {\widehat m} (y, z) f(y+x) g(z+x)dydz\bigg|^pdx\bigg)^{1/p}\\
%  &\leq &  \iint_{\R^{2n}} |{\widehat m} (y, z)|  \bigg(\int_{\R^n} |f(y+x) g(z+x)|^p\, dx \bigg)^{1/p}\, dydz \\
%  &\leq &  \iint_{\R^{2n}} |{\widehat m} (y, z)| dydz \|f\|_{p_1}\|g\|_{p_2} \\
%  &=& A_1 \|f\|_{p_1}\|g\|_{p_2}.
%\end{eqnarray*}
% \hfill{} $\Box$

\medskip
We now consider an off-diagonal case.

\begin{lemma}\label{le2.2}

\begin{itemize}
\item[(i)]    If the symbol $m(\xi, \eta)$ satisfies
\begin{eqnarray}\label{eee}
A_2=\sup_{\xi\in{\mathbb R}^n}  \Big(\int_{\R^n} | m(\xi-\eta, \eta)|^2 d\eta\Big)^{1/2} <\infty,
\end{eqnarray}
then $T_m$ maps $L^{2}(\R^n) \times L^{2}(\R^n) \to L^{2}(\R^n)$ with
$$
\|T_m(f,g)\|_2
\leq  A_2\|f\|_2 \|g\|_2.
$$

\item[(ii)] If  the symbol $m(\xi, \eta)$ is  supported on a ball
 of radius $R$,    say $B(0, R)$, and satisfies \eqref{eee},
%$$
%\sup_{\xi\in{\mathbb R}^n}  \Big(\int_{\R^n} | m(\xi-\eta, \eta)|^2 d\eta\Big)^{1/2}=A_2<\infty,
%$$
 then for all $1\leq p, q\leq 2\leq r \leq \infty$, there exists a constant $C=C_{p,q, r}$ such that
 $$
 \|T_m(f,g)\|_{r}\le C A_2 R ^{n({1\over p}+{1\over q}-{1\over r}-{1\over 2})} \|f\|_p\|g\|_q.
 $$
 \end{itemize}
\end{lemma}

\medskip

\noindent
{\bf Proof.}   The proof of (i)   follows from an application of the Plancherel identity and the Cauchy-Schwarz inequality.
\begin{eqnarray*}
\|T_m(f,g)\|_2^2
&=&  \int_{\R^n}\Big| \int_{\R^n} m(\xi-\eta, \eta) {\widehat f}(\xi-\eta) {\widehat g}(\eta) d\eta\Big|^2 d\xi\\
&\leq &  \int_{\R^n}\Big(\int_{\R^n} | m(\xi-\eta, \eta)|^2 d\eta\Big)
\Big(\int_{\R^n} |{\widehat f}(\xi-\eta) {\widehat g}(\eta)|^2 d\eta\Big) d\xi \\
&\leq &A_2^2 \|f\|_2^2 \|g\|_2^2.
\end{eqnarray*}
We now prove (ii). Since the symbol $m(\xi, \eta)$ is  supported in the ball
 $B(0, R)$, we use the Cauchy-Schwarz inequality and   Plancherel's identity  to obtain
 \begin{eqnarray*}
\|T_m(f,g)\|_{\infty}
&\leq &   \bigg\|{\mathscr F}^{-1}_\xi \bigg[\int_{\R^n}   m(\xi-\eta, \eta) {\widehat f}(\xi-\eta) {\widehat g}(\eta)d\eta \bigg]\bigg\|_{\infty}
 \\
&\lesssim &      R^{n/2}  \bigg\|\int_{\R^n}   m(\xi-\eta, \eta) {\widehat f}(\xi-\eta) {\widehat g}(\eta)d\eta\bigg\|_2\\
&\lesssim&   R^{n/2}  \|T_m(f,g)\|_2\, .
\end{eqnarray*}

In view of the support properties of $m$, in the  expression $\|T_m(f,g)\|_2$, one may
replace $ f$ and $g$ by $(\widehat f \chi_{B(0,R)})\spcheck$ and $(\widehat f \chi_{B(0,R)})\spcheck$, respectively.
 Let $r\geq 2$. It follows by interpolation and by the result  in (i) that
\begin{eqnarray*}
\|T_m(f,g)\|_r
&\lesssim &  R^{{n\over 2}(1-{2\over r})} \|T_m(f,g)\|_2
 \\
&\lesssim&   A_2 R^{{n\over 2}(1-{2\over r})} \|{\widehat f}\,  \|_{L^2(B(0,  R))} \|{\widehat g}   \|_{L^2(B(0, R))}\\
&\lesssim& A_2 R ^{n({1\over p}+{1\over q}-{1\over r}-{1\over 2})} \|f\|_p\|g\|_q.
\end{eqnarray*}
This proves (ii), and thus completes the proof of Lemma~\ref{le2.2}.
 \hfill{} $\Box$

 \medskip

 %\subsection{A useful auxiliary result}

%The next lemma will play  a crucial role in obtaining bounds for radial bilinear  multipliers.
The following  lemma is inspired by the result of Guillarmou,  Hassell, and   Sikora \cite{GHS} in the linear case.

\begin{lemma}\label{le2.6} Let $1\leq p, q \leq \infty$ and ${1/r}={1/p}+{1/q}$ and $0<r \leq \infty$.
 Suppose $T$ is a  bounded bilinear operator from $L^{p_1}({\mathbb R}^n)\times L^{q_1}({\mathbb R}^n)\to L^s({\mathbb R}^n)$
 for some $s,p_1,q_1$ satisfying $0<r\leq s$, $1\leq p_1\leq p$ and $1\leq q_1\leq q$
such that the kernel $K_T$ of $T$ satisfies
$$
\supp K_{T} \subseteq \D_{\rho}: =\big\{ (x, y, z): \, \, |x-y|<\rho, |x-z|<\rho \big\}
$$
for some $\rho>0.$ Then
there exists a constant $C=C_{r,s}>0$ such that
\begin{eqnarray}\label{e2.2-b}
\|T  \|_{L^p\times L^q\to L^r}
 \le C   \rho^{n ({1\over p_1}+{1\over q_1}-{1\over s} )}
 \| T\|_{L^{p_1}\times L^{q_1}\to L^s}.
\end{eqnarray}
\end{lemma}

\medskip

\noindent {\bf Proof.}    We fix $\rho>0$. Then  we choose a sequence of points
$(x_i)_i  $ in $ \R^n$ such that for $i\neq j$ we have
$|x_i-x_j|> \rho/10$  and $\sup_{x\in \R^n}\inf_i |x-x_i|
\le \rho/10$. Such sequence exists because $\R^n$ is separable.
Secondly, we let $B_i=B(x_i, \rho)$ and define $\widetilde{B_i}$ by the formula
$$\widetilde{B_i}=\overline{B \Big(x_i,\frac{\rho}{10}\Big)}\setminus
\bigcup_{j<i}\overline{B \Big(x_j,\frac{\rho}{10}\Big)},
$$
where $\overline{B \left(x, \rho\right)}=\{y\in \R^n \colon |x-y|
\le \rho\}$. Finally we set $\chi_i=\chi_{\widetilde B_i}$, where
$\chi_{\widetilde B_i}$ is the characteristic function of the set
${\widetilde B_i}$. Note that for $i\neq j,$
 $B(x_i, \frac{\rho}{20}) \cap B(x_j, \frac{\rho}{20})=\emptyset$. Hence
\begin{eqnarray}\label{e2.3}
 K=\sup_i\#\{j:\;|x_i-x_j|<  2\rho\} \le
  \sup_x  {|B(x, (2+\frac{1}{20})\rho)|\over
  |B(x, \frac{\rho}{20})|}=41^n<\infty.
\end{eqnarray}
It is not difficult to see that
\begin{eqnarray} \label{e2.4}
\D_{\rho}
 &\subseteq&   \bigcup_{ \substack{  i, j, k:\,  |x_i-x_j|  < 2\rho\\  \ \ \ \ \   |x_i-x_k|< 2\rho }
 } \widetilde{B}_i\times (\widetilde{B}_j \times \widetilde{B}_k)\subset \D_{4 \rho}
\end{eqnarray}
and so
$$
T(f, g) =\sum_i\sum_{\substack{  j:\,  |x_i-x_j|  < 2\rho\\  k:\,  |x_i-x_k|< 2\rho } } P_{\widetilde B_i}T (P_{\widetilde
B_j}f, P_{\widetilde B_k}g).
$$
 Let $K_r=\max\{1, K^{2(r-1)}\}$. By H\"older's inequality we have
\begin{eqnarray*}%\label{e2.5}
  \|T (f, g)\|_{r}^r
& =& \bigg\|\sum_i\sum_{\substack{  j:\,  |x_i-x_j|  < 2\rho\\
 k:\,  |x_i-x_k|< 2\rho } }P_{\widetilde B_i}T (P_{\widetilde
B_j}f, P_{\widetilde B_k}g)\bigg\|_{r}^r    \nonumber\\
& =& \sum_i \bigg\|\sum_{\substack{  j:\,  |x_i-x_j|  < 2\rho\\  k:\,  |x_i-x_k|< 2\rho } } P_{\widetilde B_i}T (P_{\widetilde
B_j}f, P_{\widetilde B_k}g)\bigg\|_{r}^r  \nonumber\\
& \lesssim &  K_r   \sum_i \sum_{\substack{  j:\,  |x_i-x_j|  < 2\rho\\  k:\,  |x_i-x_k|< 2\rho } }\| P_{\widetilde
B_i}T(P_{\widetilde
B_j}f, P_{\widetilde B_k}g)\|_{r}^r
\nonumber\\
& \lesssim &  K_r  \rho^{nr({1\over r}-{1\over s})}  \sum_i  \sum_{\substack{  j:\,  |x_i-x_j|  < 2\rho\\  k:\,  |x_i-x_k|< 2\rho } }
 \| T(P_{\widetilde
B_j}f, P_{\widetilde B_k}g)\|_{s}^r.
\end{eqnarray*}
Since $T$ is a  bounded bilinear operator from $L^{p_1}({\mathbb R}^n)\times L^{q_1}({\mathbb R}^n)\to L^s({\mathbb R}^n)$, we have
\begin{eqnarray*}
 \| T(P_{\widetilde
B_j}f, P_{\widetilde B_k}g)\|_{s} & \le  &     \|
 T\|_{L^{p_1}\times L^{q_1}\to L^s}
  \|P_{\widetilde
B_j}f\|_{p_1} \|P_{\widetilde B_k}g\|_{q_1}\nonumber\\
& \lesssim &         \rho^{n({1\over p_1}+{1\over q_1}-{1\over p}-{1\over q})}\|
 T\|_{L^{p_1}\times L^{q_1}\to L^s}  \|P_{\widetilde
B_j}f\|_{p} \|P_{\widetilde B_k}g\|_{q}\nonumber.
\end{eqnarray*}
We proceed by estimating
$$
{\mathscr E}(f,g)=\sum_i  \sum_{\substack{  j:\,  |x_i-x_j|  < 2\rho\\  k:\,  |x_i-x_k|< 2\rho } }
 \|P_{\widetilde
B_j}f\|_p^r \|P_{\widetilde B_k}g\|_{q}^r\, .
$$
Note that ${1/r}= {1/p}+{1/q}$.
  We use H\"older's inequality twice,
together with \eqref{e2.3}, to bound ${\mathscr E}(f,g)$ by
\begin{eqnarray*}
&&\hspace{-1.2cm}
  K \sum_i   \Big\{\sum_{ j:\,  |x_i-x_j|  < 2\rho  } \|P_{\widetilde B_j}f\|_p^{p} \Big\}^{r/p}
 \Big\{ \sum_{ k:\,  |x_i-x_k|  < 2\rho  } \|P_{\widetilde
 B_k}g\|_q^{q} \Big\}^{r/q} \\
 &\leq &
  K  \Big\{\sum_i    \sum_{ j:\,  |x_i-x_j|  < 2\rho  } \|P_{\widetilde
B_j}f\|_p^{p} \Big\}^{r/p}
\Big\{ \sum_i  \sum_{ k:\,  |x_i-x_k|  < 2\rho  } \|P_{\widetilde
B_k}g\|_q^{q} \Big\}^{r/ q} \\
&\leq &   K^2\Big\{ \sum_j   \|P_{\widetilde
B_j}f\|_p^p \Big\}^{r/p}
\Big\{ \sum_k  \|P_{\widetilde
B_k}g\|_q^q \Big\}^{r/q} \\
&\leq &  K^2\|f\|_p^r
 \| g\|_q^r.
\end{eqnarray*}
This estimate  combined with the previously obtained estimate for  $\|T(f,g)\|_r^r$ in terms of
${\mathscr E}(f,g)$ yields  \eqref{e2.2-b}.  The proof   is now complete. \hfill{}
$\Box$

 \bigskip

It will be useful to apply   Lemma \ref{le2.6} for operators, which do not have   such perfect localization properties.
For such, we have the following version:

\begin{lemma}\label{le2.6-bis} Let $1\leq p, q \leq \infty$ and ${1/r}={1/p}+{1/q}$ and $0<r \leq \infty$.
 Suppose $T$ is a  bounded bilinear operator from $L^{p_1}({\mathbb R}^n)\times L^{q_1}({\mathbb R}^n)\to L^s({\mathbb R}^n)$
 for some $s,p_1,q_1$ satisfying $0<r\leq s$, $1\leq p_1\leq p$ and $1\leq q_1\leq q$
such that the kernel $K_T$ of $T$ satisfies
$$
\left|K_T(x,y,z) \right| \lesssim \rho^{-d} \left(1+ \rho^{-1}|x-y|+\rho^{-1}|x-z|\right)^{-M}
$$
for some $\rho\geq 1$, $d>0$ and every large enough integer $M>0$. Then for every $\epsilon>0$ (as small as we want) and $N>0$ (as large as we want)
there exists a constant $C=C_{r,s,\epsilon}>0$ such that
\begin{eqnarray}\label{e2.2-bb}
\|T  \|_{L^p\times L^q\to L^r}
 \le C   \rho^{ \epsilon + n ({1\over p_1}+{1\over q_1}-{1\over s} )}
 \| T\|_{L^{p_1}\times L^{q_1}\to L^s} + \rho^{-N}.
\end{eqnarray}
\end{lemma}

\medskip

\noindent {\bf Proof.} The proof is very similar to the previous one.
Let us fix $\epsilon>0$ and consider a collection of points $(x_i)_i$ in $\R^n$ such that for $i\neq j$ we have
$|x_i-x_j|> \rho^{1+\epsilon}/10$  and $\sup_{x\in \R^n}\inf_i |x-x_i|
\le \rho^{1+\epsilon} /10$.
Then, with the previous notation, we have
\begin{eqnarray*}
&&\hspace{-0.6cm}\|T (f, g)\|_{r}  \\
&& \lesssim   \bigg\|\sum_i\sum_{\substack{  j:\,  |x_i-x_j|  < 2\rho^{1+\epsilon}\\
 k:\,  |x_i-x_k|< 2\rho^{1+\epsilon} } }P_{\widetilde B_i}T (P_{\widetilde
B_j}f, P_{\widetilde B_k}g)\bigg\|_{r}  + \bigg\|\sum_i\sum_{\substack{  j,k:\,  |x_i-x_j|  > 2\rho^{1+\epsilon}\\
 \textrm{or }  |x_i-x_k|> 2\rho^{1+\epsilon} } }P_{\widetilde B_i}T (P_{\widetilde
B_j}f, P_{\widetilde B_k}g)\bigg\|_{r} \\
&&:=   I+II.
\end{eqnarray*}
For the  $I$, we  repeat exactly the same reasoning as for Lemma \ref{le2.6} (since it
corresponds to the diagonal part), by replacing $\rho$ by $\rho^{1+\epsilon}$. So we obtain
$$ I \lesssim \rho^{(1+\epsilon) n ({1\over p_1}+{1\over q_1}-{1\over s} )}
 \| T\|_{L^{p_1}\times L^{q_1}\to L^s},$$
which is as claimed since $\epsilon$ can be chosen as small as we want.

We now deal with the second quantity $II$. We have
\begin{align*}
 II & \leq \sum_{2^\ell \geq \rho^\epsilon} \bigg\|\sum_i \sum_{\substack{  j,k:\, \\ |x_i-x_j| + |x_i-x_k| \simeq \rho 2^\ell } }
 P_{\widetilde B_i}T (P_{\widetilde
B_j}f, P_{\widetilde B_k}g) \bigg\|_{r} \\
 & \lesssim \sum_{2^\ell \geq \rho^\epsilon} \rho^{-d} 2^{-\ell M} (\rho 2^\ell)^{3n} \ \Big\| \sum_{j}   |P_{\widetilde
B_j}f  | \Big\|_{p} \Big\| \sum_{k}  |P_{\widetilde B_k}g  | \Big\|_{q},
\end{align*}
where we used that for $j,k$ fixed, there is at most $(\rho 2^\ell)^{n}$ points $x_i$ satisfying
$$
|x_i-x_j| + |x_i-x_k|  \simeq \rho 2^\ell
$$
and the pointwise estimate of the bilinear kernel.
So we conclude that
\begin{align*}
 II & \lesssim \|f\|_{p} \|g\|_q \bigg( \sum_{2^\ell \geq \rho^\epsilon} \rho^{-d} 2^{-\ell M} (\rho 2^\ell)^{3n} \bigg) \\
 & \lesssim \rho^{-d-\epsilon M +3n(1+\epsilon)}\|f\|_{p} \|g\|_q,
\end{align*}
which is also as claimed since $M$ can be chosen as large as we want. \hfill{}
$\Box$

\medskip

 \section{Compactly supported bilinear multipliers}
\setcounter{equation}{0}

 In this  section  we assume   $n\geq 2 $ and we are concerned with the
 boundedness of compactly supported bilinear multipliers. We focus attention to radial such multipliers. These can be written in the form
\begin{eqnarray}\label{e3.1}
  \iint_{\R^{2n}}  e^{2\pi ix\cdot(\xi+\eta)} m_0(|\xi|, |\eta|) {\widehat f}(\xi) {\widehat g}(\eta)d\xi d\eta
\end{eqnarray}
for $f, g\in{\mathscr S}(\R^n)$,
where $m_0\in L^{\infty}(\R^+\cup\{0\}\times \R^+\cup\{0\})$. This is exactly the bilinear multiplier
operator $T_m(f,g)$, where $m_0(|\xi|, |\eta|)=  m(\xi_1, \dots , \xi_n, \eta_1, \dots, \eta_n)  $,
$\xi=(\xi_1, \dots , \xi_n)$ and $\eta=(\eta_1, \dots, \eta_n)$.
%With some abuse of notation, we will denote $T_{m}$ by $T_{m_0}$ without causing any confusion.

\begin{lemma} \label{lem:band} Let $m_0$ be an even function on $\R^2$ whose  Fourier transform $\widehat{m_0}$  is supported in $[-L,L]^2$.
Then the Fourier transform of the biradial
 function $m(x,y):=m_0(|x|,|y|)$ on $\R^{2n}$ is    supported in $[-L,L]^{2n}$.
\footnote{We give here a proof, using the finite speed propagation property of the wave propagator. Actually in the linear
 framework, the claim can be rephrased as follows: a Fourier band-limited function is also a Hankel band-limited function,
 for the ``$J_0$'' Hankel transform. We also refer the reader to \cite{Rawn,C} for another approach  to this question using the Hankel transform.}
\end{lemma}

\medskip

\noindent
{\bf Proof.} We have
\begin{align*}
 \widehat{m}(\xi,\eta) & = \iint_{\R^{2n}} e^{2i\pi (x \cdot\xi + y\cdot \eta)} m_0(|x|,|y|) dx dy \\
 & = \iint_{[0,\infty)^2} m_0(u,v) R_u(\xi) R_v(\eta) du\, dv \\
 & =  m_0(\sqrt{-\Delta},\sqrt{-\Delta}) (\delta_0,\delta_0) (\xi,\eta) \\
 & = K_{m_0(\sqrt{-\Delta},\sqrt{-\Delta})}(\xi,\eta),
\end{align*}
where $K_{m_0(\sqrt{-\Delta},\sqrt{-\Delta})}$ is the bilinear kernel of the bilinear operator
 $m_0(\sqrt{-\Delta},\sqrt{-\Delta})$. Expressing $m_0$ in terms of its 2-dimensional Fourier transform yields
\begin{align*}
K_{m_0(\sqrt{-\Delta},\sqrt{-\Delta})}(\xi,\eta) & = \int_{\R^2} \widehat{m_0}(s,t)
K_{e^{2i\pi s \sqrt{-\Delta}} e^{2i\pi t\sqrt{-\Delta}}} (\xi,\eta) du\, dv \\
 & = \int_{[-L,L]^2} \widehat{m_0}(s,t) K_{e^{2i\pi s \sqrt{-\Delta}}}(\xi)
 K_{e^{2i\pi t\sqrt{-\Delta}}}(\eta) du\, dv.
\end{align*}
Then using finite speed propagation property of the wave propagator, we know that for
every $s\in\R$, the kernel $K_{e^{2i\pi s \sqrt{-\Delta}}}$ is supported on $[-|s|,|s|]^n$.
Hence, we conclude that $K_{m_0(\sqrt{-\Delta},\sqrt{-\Delta})}$ is supported on $[-L,L]^{2n}$.
\hfill{} $\Box$

 \medskip

 \subsection{Bilinear restriction-extension operators}

   For $f \in{\mathscr S}(\R^n)$ recall the restriction-extension operator
\begin{eqnarray}\label{e3.2}\hspace{1cm}
{\mathscr R}_{\lambda}f(x)=  {\lambda^{n-1}  }\int_{{\mathbb S}^{n-1}}  e^{2\pi i\lambda x\cdot \omega }
 {\widehat f}(\lambda \omega ) d\omega
 %= (2\pi)^{-n}\int_{\lambda {\bf S}^{n-1}} e^{ix \omega} \widehat{f}(\omega) d\omega
 \ \  \ \ \lambda>0
 \end{eqnarray}
  in the linear setting.
 In the sequel  we set
$$
 {\mathscr R}_{\lambda_1, \lambda_2}(f, g)(x)=   {\mathscr R}_{\lambda_1}f(x){\mathscr R}_{\lambda_2}g(x).
$$

\begin{lemma} \label{le3.1} Let $ m(\xi,\eta) :=m_0(|\xi|, |\eta|)  $.
For $f, g\in{\mathscr S}(\R^n),$
we have the following formula:
\begin{eqnarray*}
T_m(f, g)(x) =   \int_0^{\infty} \int_0^{\infty} m_0(\lambda_1, \lambda_2)
 {\mathscr R}_{\lambda_1, \lambda_2}(f, g)(x) d\lambda_1 d\lambda_2.
\end{eqnarray*}
\end{lemma}

\medskip

\noindent
{\bf Proof.} The proof can be obtained by expressing $ T_m(f, g)(x)$ in polar coordinates.
 \hfill{} $\Box$

\medskip

To study boundedness of the bilinear restriction-extension operator ${\mathscr R}_{\lambda_1, \lambda_2},$
we first recall some properties of the operator ${\mathscr R}_{1}$ in the linear setting.
Let $d\sigma$ denote  surface measure on the unit sphere ${\mathbb S}^{n-1}$.
In view of  the theory of Bessel function
 (see   page  428  of \cite{Gra}),
 \begin{eqnarray}\label{e3.3}
{\mathscr R}_1f(x)=   \int_{{\mathbb S}^{n-1}} e^{2\pi ix\cdot \theta} {\widehat f}(\theta) d\theta= {\widehat d\sigma}  \ast f (x),
% = c {\widehat{ \big({\breve f}d\theta \big)}}(x),
\end{eqnarray}
where
\begin{eqnarray} \label{e3.4}
{\widehat d\sigma} (x) =\int_{{\bf S}^{n-1}} e^{2\pi ix \cdot \omega} d \omega =    {2\pi\over |x|^{n-2\over 2}} J_{n-2\over 2}(2\pi |x|)
\end{eqnarray}
and
\begin{eqnarray}\label{e3.5}
 J_\zeta(t)&=& \frac{2}{\Gamma(1/2)}\frac{(t/2)^\zeta}{\Gamma(\zeta+1/2)}\int_0^1(1-u^2)^{\zeta-1/2}\cos
(ut)du,\ \ \ {\rm Re} (\zeta) >-{1\over 2}.
\end{eqnarray}

The problem of $L^p$-$L^q$ boundedness of ${\mathscr R}_1$ has been studied by several authors (see for instance,
\cite{Bk}, \cite {Br} and  \cite{Gu}).
The first results in this
direction were obtained by Tomas and  Stein \cite{Tom}, \cite{St2}; they showed that ${\mathscr R}_1$ is bounded from
$L^p(\R^n)$ to $L^{p'}(\R^n)$ for $p=(2n+2)/(n+3)$ and $p'=p/(p-1)$, which implies the
sharp $L^p-L^2$ restriction theorem for the sphere ${\mathbb S}^{n-1}.$ To describe all pairs $(p,q)$
such that the operator ${\mathscr R}_1$ on $\R^n$ is bounded from
$L^p(\R^n)$ to $L^q(\R^n)$,
%Background related to the subject can be founded in \cite {Br}, \cite{Sog3},   \cite{BMO} and \cite{Bk}.
we define vertices in the square $[0,1]\times [0,1]$ by setting
\begin{eqnarray*}
A(n)=\Big({n+1\over 2n}, 0\Big), \ \ \ \ \ \  B(n)= \Big({n+1\over 2n}, {n-1\over 2n}-{n-1\over n^2+n}\Big),\nonumber\\[4pt]
A'(n)=\Big(1, {n-1\over 2n}\Big), \ \ \ \ \ \  B'(n)= \Big({n+1\over 2n}+{n-1\over n^2+n},  {n-1\over 2n}\Big).
\end{eqnarray*}
Let $\Delta(n)$ be the closed pentagon with vertices $A(n), B(n), B'(n), A'(n)$, $(1,0)$
from which closed line segments $[A(n), B(n)], [A'(n), B'(n)]$ are removed. Namely,
\begin{equation}
\label{e3.6}
\Delta(n)=\left\{
\begin{array} {ll}
\Big(\dfrac{1}{p}, \dfrac{1}{q}\Big)\in [0,1]\times [0,1]:
 & 0 \leq \dfrac{1}{q}\leq \dfrac{1}{p}\leq 1,\ \  \dfrac{1}{p}-\dfrac{1}{q}\geq \dfrac{2}{n+1}  ,
  \\ [13pt]
   & \ \ \ \
  \dfrac{1}{p}>\dfrac{n+1}{2n} {\over }, \ \ \dfrac{1}{q}< \dfrac {n-1}{2n}
\end{array}
  \right\}
\end{equation}

\begin{proposition} \label{prop3.2}
Let ${\mathscr R}_1$ be defined as in \eqref{e3.2}. There exists a constant $C=C_{p,q}>0$, independent of $f$,
such that
$$
\|{\mathscr R}_1(f)\|_q\leq C\|f\|_p
$$
if and only if  $(1/p, 1/q)$ in $\Delta(n)$.
\end{proposition}

\medskip

 \noindent
{\bf Proof.} For the proof of Proposition~\ref{prop3.2}, we refer it to Remark 1,  p. 497,   \cite{Gu}. See also \cite{Bk},
and \cite{Br}.
 \hfill{} $\Box$

\medskip

As a consequence of Proposition~\ref{prop3.2}, we obtain the following result.

\begin{proposition}  \label{prop3.3}
(i) Let   $ {1/s}={1/q_1}+ {1/q_2}$,  $0<s\leq\infty$ and let $  (1/p_1, 1/q_1)$  and $(1/p_2, 1/q_2)$ be both in
$\Delta(n)$ as defined as in \eqref{e3.6}.
   For every $\lambda_1, \lambda_2>0$, the bilinear  restriction-extension operator
     ${\mathscr R}_{\lambda_1, \lambda_2}$ is bounded from $L^{p_1}(\R^n)\times L^{p_2}(\R^n)$ to $L^s(\R^n)$ such that
\begin{eqnarray*}
 \|{\mathscr R}_{\lambda_1, \lambda_2}(f, g)\|_s\lesssim
 \lambda_1^{n({1\over p_1}-{1\over q_1})-1} \lambda_2^{n({1\over p_2}-{1\over q_2})-1}\|f\|_{p_1} \|g\|_{p_2}.
 \end{eqnarray*}

\noindent (ii) In the  endpoint case $s=2$ and $p_1=p_2=1$ we have    ,
\begin{equation}
\|{\mathscr R}_{\lambda_1, \lambda_2}(f, g)\|_2 \lesssim  \lambda_1^{n-\frac{3}{2}} \lambda_2^{\frac{n-1}{2}} \| f \|_{1} \|g\|_{1}, \label{eq:endpoint}
\end{equation}
assuming $\lambda_1<<\lambda_2$. This
  corresponds to the   result in (i) with $q_2= \frac{2n}{n-1}$.
\end{proposition}

 \medskip

\noindent
{\bf Proof.}
(i)  By H\"older's inequality,
\begin{eqnarray}\label{e3.7}
 \|{\mathscr R}_{\lambda_1, \lambda_2}(f, g)\|_s= \|{\mathscr R}_{\lambda_1} f  {\mathscr R}_{\lambda_2} g\|_s
 & \leq& \|{\mathscr R}_{\lambda_1} f   \|_{q_1}  \|  {\mathscr R}_{\lambda_2} g\|_{q_2},
\end{eqnarray}
\noindent where ${1/s}={1/q_1}+ {1/q_2}.$
Note that  by Proposition~\ref{prop3.2}, we  obtain that if  both $  (1/p_1, 1/q_1)$  and $(1/p_2, 1/q_2)$ are in $\Delta(n)$, then
\begin{eqnarray}\label{e3.8}
 \|{\mathscr R}_{\lambda_1} f   \|_{q_1}
 \lesssim  \lambda_1^{n({1\over p_1}-{1\over q_1})-1} \|f\|_{p_1}
  \end{eqnarray}
  and
\begin{eqnarray}\label{e3.9}
 \|{\mathscr R}_{\lambda_2} g   \|_{q_2}
 \lesssim \lambda_2^{n({1\over p_2}-{1\over q_2})-1} \|g\|_{p_2}.
  \end{eqnarray}
The desired estimate   now follows from    \eqref{e3.7},  \eqref{e3.8} and \eqref{e3.9}.

(ii) We now prove (\ref{eq:endpoint}).
Let $\mathbb B$ be the unit ball in $\R^n$.
First, ${\mathscr R}_{\lambda_1, \lambda_2}(f, g)$ has a spectrum included in
$$
Sp:=\lambda_1 {\mathbb S}^{n-1} + \lambda_2 {\mathbb S}^{n-1} = (\lambda_1+\lambda_2){\mathbb B} \setminus (\lambda_2-\lambda_1){\mathbb B}
$$
which has a $n$-dimensional measure
\begin{equation}
 |Sp|=|\lambda_1 {\mathbb S}^{n-1} + \lambda_2 {\mathbb S}^{n-1}| \lesssim \lambda_1 \lambda_2^{n-1}
\label{eq:sp} \end{equation}
since $\lambda_1<<\lambda_2$.
Moreover, since $f,g\in L^1$ then by Plancherel equality, we have
 \begin{align*}
   \|{\mathscr R}_{\lambda_1, \lambda_2}(f, g)\|_2 & = \|{\mathscr F} \left[{\mathscr R}_{\lambda_1, \lambda_2}(f, g)\right] \|_2 \\
 & = \| \widehat{{\mathscr R}_{\lambda_1}(f)} \ast \widehat{ {\mathscr R}_{\lambda_2}(g)} \|_2 \\
 & = \sup_{\genfrac{}{}{0pt}{}{h\in L^2}{\|h\|_2=1}} \iint_{\lambda_1 {\mathbb S}^{n-1} \times \lambda_2
 {\mathbb S}^{n-1}} \widehat{{\mathscr R}_{\lambda_1}(f)} (\xi) \widehat{ {\mathscr R}_{\lambda_2}(g)} (\eta) h(\xi+\eta) d\xi d\eta \\
 & \leq \|f\|_1 \|g\|_1 \sup_{\genfrac{}{}{0pt}{}{h\in L^2}{\|h\|_2=1}} \iint_{\lambda_1 {\mathbb S}^{n-1} \times
  \lambda_2 {\mathbb S}^{n-1}} |h(\xi+\eta)| d\xi d\eta \\
 & \leq \|f\|_1 \|g\|_1 \lambda_1^{n-2} \sup_{\genfrac{}{}{0pt}{}{h\in L^2}{\|h\|_2=1}} \int_{Sp} |h(\omega)| d\omega.
  \end{align*}
At the last inequality, we used that for every $\omega\in Sp$
$$ {\textrm{meas}}\left(\{(\xi,\eta) \in \lambda_1 {\mathbb S}^{n-1}
 \times \lambda_2 {\mathbb S}^{n-1}:\,\, \xi+\eta=\omega\}\right) \lesssim \lambda_1^{n-2},$$
where ${\textrm{meas}}$ denotes  $(n-2)$-dimensional Hausdorff measure.
Finally,  via (\ref{eq:sp}) we obtain
 \begin{align*}
   \|{\mathscr R}_{\lambda_1, \lambda_2}(f, g)\|_2 & \lesssim \|f\|_1 \|g\|_1 \lambda_1^{n-2} |Sp|^{\frac{1}{2}}
   \lesssim \|f\|_1 \|g\|_1 \lambda_1^{n-2} (\lambda_1 \lambda_2^{n-1})^{\frac{1}{2}},
  \end{align*}
which yields (\ref{eq:endpoint}).
 \hfill{} $\Box$

 \medskip

 \subsection{Restriction-extension estimates imply  bilinear multiplier estimates}

For every $n\geq 2$, set
$$
a_n={n+1\over 2n}, \ \ \ {\rm and}\ \ \ b_n={n+1\over 2n} +{n-1\over n^2+n}.
$$
Let  $\epsilon>0$ and for every  $1\leq p_1, p_2 < {2n\over n+1}$ and   $ {1\over p}={1\over p_1}+ {1\over p_2}$,
 define
 \begin{eqnarray} \label{e3.10} \hspace{2cm}
\alpha(p_1, p_2, \epsilon)=
\left\{
\begin{array}{lll}
{4\over n+1}, &\ \ {\rm if}\ \ \big({1\over p_1}, {1\over p_2}\big)\in \big(a_n, b_n\big)\times
 \big(a_n, b_n\big) \\[10pt]
{2\over n+1} -{n-1\over 2n}+{1\over p_2} +\epsilon, &\ \ {\rm if}\ \  \big({1\over p_1}, {1\over p_2}\big)
\in \big(a_n, b_n\big)
\times \big[b_n, 1\big]  \\[10pt]
{2\over n+1} -{n-1\over 2n}+{1\over p_1}+\epsilon,&\ \ {\rm if}\ \ \big({1\over p_1}, {1\over p_2}\big)\in
 \big[ b_n, 1\big]\times
  \big(a_n, b_n\big) \\[10pt]
{1\over p_1}+{1\over p_2}-{n-1\over  n}+\epsilon, &\ \ {\rm if}\ \ \big({1\over p_1}, {1\over p_2}\big)\in
\big[ b_n, 1\big]\times
  \big[ b_n, 1\big].
 \end{array}
 \right.
 \end{eqnarray}
For simplicity, we will write $\alpha(p_1, p_2)$ instead of $\alpha(p_1, p_2, 0).$
   \medskip

Now we prove the following  result.

 \begin{theorem}\label{th3.4} Let  $1\leq p_1, p_2 < 2n/(n+1)$ and
  $ {1/p}={1/p_1}+ {1/p_2}$.
Suppose  that
  $m_0$ is an even bounded function supported in $ [-1,1]\times [-1, 1]$ that lies in $ W^{\beta,\, 1}(\R^2)$,
  for some $\beta>n\alpha(p_1, p_2) $.
Let  $m(\xi,\eta):=  m_0(|\xi|,|\eta|)$.
Then $T_{m}$ is a bounded operator from  $L^{p_1}(\R^n)\times L^{p_2}(\R^n)$ to $L^p(\R^n)$ and we
have
\begin{equation}
\label{e3.11}  \|T_{ m}\|_{L^{p_1}\times L^{p_2}\to L^p} \leq
C\|m_0\|_{W^{\beta, 1}({\R^2})}.
\end{equation}
\end{theorem}

\medskip

\noindent
{\bf Proof.}
Let $\phi\in C_c^{\infty}({\R})$ be   an even   function with  $\supp \phi\subseteq \{ t: \ 1/4\leq |t|\leq 1\}$ and
 %such that $\phi(\lambda)=1$ for all $\lambda\in [{3\over 8}, {7\over 8}]$, and
$$
\sum_{\ell\in \ZZ} \phi(2^{-\ell} \lambda)=1 \ \  \ \ \forall
{\lambda>0}.
$$
Then we  set $\phi_0(\lambda)= 1-\sum_{\ell\geq 1} \phi(2^{-\ell} \lambda)$,
\begin{eqnarray}\label{e3.12}\hspace{0.5cm}
m^{(0)}_0(\lambda_1, \lambda_2) =
\iint_{\R^2}
 \phi_0(\sqrt{ |t_1|^2 +|t_2| ^2}\,) \widehat{m_0}(t_1, t_2)
e^{2\pi i (t_1\lambda_1+t_2\lambda_2) }
  \;dt_1 dt_2
\end{eqnarray}
and
\begin{eqnarray}\label{e3.13}\hspace{0.5cm}
m^{(\ell)}_0(\lambda_1, \lambda_2) = \iint_{\R^2}
\phi\big(2^{-\ell} \sqrt{ |t_1|^2 +|t_2| ^2}\,\big) \widehat{m_0}(t_1, t_2)
e^{2\pi i (t_1\lambda_1+t_2\lambda_2) }
\;dt_1 dt_2.
\end{eqnarray}
Note that in view of Fourier inversion and of the preceding decomposition we have that
\begin{eqnarray}\label{e3.14}
m_0(\lambda_1, \lambda_2)&=&\sum_{\ell\geq 0}  m^{(\ell)}_0(\lambda_1, \lambda_2).
\end{eqnarray}
Consequently,
\begin{eqnarray}\label{e3.15}
 T_{m}(f,g)(x)
&=&  \iint_{\R^{2n}}  e^{2\pi ix\cdot (\xi+\eta)}
 m_0 (|\xi|, |\eta|) {\widehat f}(\xi) {\widehat g}(\eta)d\xi d\eta \nonumber\\
 &=&\sum_{\ell\geq 0}  T_{m^{(\ell)} }(f,g)(x),
 \end{eqnarray}
 where $m^{(\ell)}(\xi,\eta) = m_0^{(\ell)}(|\xi|,|\eta|)$ for $\ell\ge 0$ and
 \begin{eqnarray}\label{e3.16}
  T_{m^{(\ell)}}(f,g)(x)=  \iint_{\R^{2n}}  e^{2\pi ix\cdot (\xi+\eta)}
  m^{(\ell)}_0 (|\xi|, |\eta|) {\widehat f}(\xi) {\widehat g}(\eta)d\xi d\eta.
\end{eqnarray}
It follows from the support properties of $\phi$ and Lemma \ref{lem:band} that the kernel of $T_{ m^{(\ell)}} $ is supported in
$$
  \D_{2^{\ell} }=\big\{ (x, y, z):\, \,
 |x-y|<2^{\ell}  , |x-z|<2^{\ell}  \big\}.
$$

   Recall the set  $\Delta(n)$    given in \eqref{e3.6}. We observe that
if $\epsilon>0$ is small enough, then there exist $({1/p_1}, 1/q_1) \in \Delta(n), ({1/p_2}, 1/q_2) \in \Delta(n) $
such that
$$
{1\over p_1} + {1\over p_2} -{1\over q_1}-{1\over q_2}=\alpha(p_1, p_2, \epsilon).
$$
Note that $p<1.$ Let $ {1\over q_1}+{1\over q_2}={1\over s}$ and   so    $s>1 $.
Then Lemma~\ref{le2.6} yields the existence of a constant $C=C_{p,s}$ such that
\begin{eqnarray} \label{e3.17}\hspace{1cm}
\|T_{ m^{(\ell)}}\|_{L^{p_1}\times L^{p_2}\to L^p}
 &\le& C\,  2^{\ell  n   ({1\over p}-{1\over s})}
 \| T_{ m^{(\ell)}}  \|_{L^{p_1}\times L^{p_2}\to L^s},
\end{eqnarray}
which yields
\begin{eqnarray}\label{e3.19}
 \| T_{  m } \|_{L^{p_1}\times L^{p_2}\to L^p}&\leq& \Big\| \sum\limits_{\ell\geq 0}
 T_{   m^{(\ell)}  } \Big\|_{L^{p_1}\times L^{p_2}\to L^p}
 \nonumber\\
  &\lesssim &
\sum\limits_{\ell\geq 0}   2^{\ell \theta + \ell n   ({1\over p}-{1\over s})}
  \| T_{  m^{(\ell)}  } \|_{L^{p_1}\times L^{p_2}\to L^s}
 \end{eqnarray}
for some constant  $\theta\in (0, (\beta-n\alpha(p_1, p_2))/2)$.

Since $m^{(\ell)}$ is not compactly supported we choose a smooth function
 $\psi$ supported in $(-8, 8)$ such that $\psi(\lambda)=1$ for $\lambda \in (-4,4)$. We set
 $\Psi(x_1,x_2) = \psi(|x_1|+|x_2|)$ for $x_1,x_2\in \R^n$
and we note that
\begin{eqnarray} \label{e3.20} \hspace{1cm}
\| T_{ m^{(\ell)}}  \|_{L^{p_1}\times L^{p_2}\to L^s}\leq \| T_{  \Psi  m^{(\ell)}  } \|_{L^{p_1}\times L^{p_2}\to L^s}
+\| T_{ (1-\Psi)   m^{(\ell)}  }  \|_{L^{p_1}\times L^{p_2}\to L^s},
\end{eqnarray}
where $(\Psi  m^{(\ell)})(\xi, \eta)=:\psi(|\xi|+|\eta|)m^{(\ell)}_0(|\xi|, |\eta|).$

To estimate $ \| T_{  \Psi  m^{(\ell)}  } \|_{L^{p_1}\times L^{p_2}\to L^s},$
 we apply Lemma~\ref{le3.1},   together with
 Minkowski's inequality ($s>1$),  and Proposition~\ref{prop3.3}    to
obtain
\begin{eqnarray}\label{e1111}
\| T_{  \Psi  m^{(\ell)}  } \|_{L^{p_1}\times L^{p_2}\to L^s}
& \leq &
  \int_0^{\infty} \int_0^{\infty}   |  \psi (\lambda_1+\lambda_2)  m^{(\ell)}_0  (\lambda_1, \lambda_2)|
  \| {\mathscr R}_{\lambda_1, \lambda_2}  \|_{L^{p_1}\times L^{p_2}\to L^s}
    d\lambda_1 d\lambda_2\nonumber  \\
&  \lesssim &
 \int_0^{8 } \int_0^{8 }   |   m^{(\ell)}_0  ( \lambda_1, \lambda_2)|
 \lambda_1^{n({1\over p_1}-{1\over q_1})-1} \lambda_2^{n({1\over p_2}-{1\over q_2})-1}   d\lambda_1 d\lambda_2\nonumber\\
&  \lesssim &
 \,   \|     m^{(\ell)}  \|_1
%&  \leq & C \,  \|  \psi  m^{(\ell)} \|_1.
\end{eqnarray}
since $\frac{1}{p_i}-\frac{1}{q_i} \ge \frac{2}{n+1} $, thus $n(\frac{1}{p_i}-\frac{1}{q_i})-1 \ge \frac{n-1}{n+1}  \ge 0$
in view of the fact that $(\frac{1}{p_i},\frac{1}{q_i}) \in \Delta(n)$.

Notice that ${1\over p}-{1\over s}=
{1\over p_1} + {1\over p_2} -{1\over q_1}-{1\over q_2}=\alpha(p_1, p_2, \epsilon) $.
 Then we have
\begin{eqnarray}\label{e3.21}
\sum_{\ell\geq 0}  2^{\ell \theta + \ell n   ({1\over p}-{1\over s})}
  \| T_{  \Psi  m^{(\ell)}_0 } \|_{L^{p_1}\times L^{p_2}\to L^s}
  &\lesssim&    \sum_{\ell\geq 0}  2^{\ell(n\alpha(p_1, p_2, \epsilon)+\theta)}  \|    m^{(\ell)}_0  \|_1  \nonumber\\
  &\lesssim &    \|m_0\|_{B^{n\alpha(p_1, p_2, \epsilon)+\theta}_{1,1}},
\end{eqnarray}
where   the last inequality follows from the definition
Besov space.
See, e.g., \cite[Chap.~VI ]{BL}.
Recall also that if $\beta >n\alpha(p_1, p_2, \epsilon) +\theta$ then
$$
W^{\beta,1}({\R^2})\subseteq
B_{1, \, 1}^{n\alpha(p_1, p_2, \epsilon)+\theta}({\R^2})
$$
 and $\|m_0\|_{B_{1,\, 1}^{n\alpha(p_1, p_2, \epsilon)+\theta}({\R^2})}\le C_\beta
\|m_0\|_{W^{\beta,1}({\R^2})}$, see again
\cite{BL}.

Next we obtain bounds for
   $\| T_{ (1-\Psi)   m^{(\ell)} }  \|_{L^{p_1}\times L^{p_2}\to L^s}$.
Since the function $1-\psi$ is supported outside the
interval $(-4,4)$,
we can choose  a function $\eta\in C_c^{\infty}(4, 16)$  such that
$$
1=\psi(t)+\sum_{k\geq 0}\eta(2^{-k}t)
$$
 for all $t>0$.
Hence for $\lambda_1,\lambda_2>0$ we have
\begin{eqnarray*}
\big( 1-\psi  ( \lambda_1+\lambda_2) \big)m^{(\ell)}_0(\lambda_1, \lambda_2)= \sum_{k\geq 0}
\eta ( 2^{-k}(\lambda_1 + \lambda_2) )m^{(\ell)}_0 (\lambda_1, \lambda_2).
\end{eqnarray*}
We then apply an argument as in \eqref{e1111} to show that
\begin{eqnarray} \label{e222}\hspace{1cm}
&&\hspace{-1cm}  \| T_{  (1-\Psi)   m^{(\ell)} }   \|_{L^{p_1}\times L^{p_2}\to L^s} \nonumber\\
  &&\lesssim
    \sum\limits_{k\geq 0}  \int_0^{2^{k+4}}\!\!\!\int_0^{2^{k+4}}   | \eta  (2^{-k}  ( \lambda_1+\lambda_2)) m^{(\ell)}_0 (\lambda_1, \lambda_2)|
 \lambda_1^{n({1\over p_1}-{1\over q_1})-1} \lambda_2^{n({1\over p_2}-{1\over q_2})-1}  d\lambda_1 d\lambda_2.
\end{eqnarray}
Observe that
\begin{eqnarray*}
 m^{(\ell)}_0 ( \lambda_1,  \lambda_2)
= \iint\limits_{[-1,1]^2} m_0(s_1,s_2) \bigg[ \iint\limits_{\R^2} \phi(2^{-\ell}\sqrt{|t_1|^2+|t_2|^2}\, )
e^{2\pi i   ( (\l_1-s_1)t_1+(\l_2-s_2)t_2)}dt_1dt_2 \bigg] ds_1ds_2.
\end{eqnarray*}
%Since   $ \supp \eta\subseteq [4, 16], \supp\phi\subseteq [1/4, 1] $  and  $\supp m_0\subseteq  [-1, 1]\times [-1, 1]$,
We can integrate by parts $M$ times to obtain
\begin{eqnarray*}
    | \eta (2^{-k} ( \lambda_1+\lambda_2))   m^{(\ell)}_0 ( \lambda_1,  \lambda_2)|
 \le  C_M2^{-(\ell+k)M + 2\ell}\|m_0\|_1.
\end{eqnarray*}
Substituting this  back into    \eqref{e222} with
$M$ sufficiently large  such that
$$
n\Big({1\over p_1}+{1\over p_2}-{1\over q_1}-{1\over q_2}\Big)-M +2 +\theta<0,
$$
 we obtain
\begin{eqnarray*}
\| T_{ (1-\Psi)   m^{(\ell)}  }  \|_{L^{p_1}\times L^{p_2}\to L^s}
  &\leq&   C_M 2^{-\ell (M-2)}\|m_0\|_1
    \sum_{k\geq 0}   2^{-kM+ kn({1\over p_1}+{1\over p_2}-{1\over q_1}-{1\over q_2})}  \\
	&\leq&   C_M 2^{-\ell (M-2)}\|m_0\|_1,
	\end{eqnarray*}
which yields
\begin{eqnarray}\label{e3.22}
 \sum_{\ell\geq 0}2^{\ell \theta + \ell n   ({1\over p}-{1\over s})}
\| T_{  (1-\Psi)   m^{(\ell)}  }  \|_{L^{p_1}\times L^{p_2}\to L^s}
    &\lesssim&     \|m_0\|_1.
\end{eqnarray}

 Finally,  \eqref{e3.11} follows from  \eqref{e3.15},     \eqref{e3.19},  \eqref{e3.20}, \eqref{e3.21} and     \eqref{e3.22}.
 This completes  the proof of Theorem~\ref{th3.4}.
  \hfill{} $\Box$

 \medskip

\begin{remark}
The previous proof relies on a bilinear spherical decomposition of the symbol. The bilinear restriction operator
 ${\mathscr R}_{\lambda_1, \lambda_2}$ is not well-defined on $L^2(\R^n)\times L^2(\R^n)$ and so one cannot use these elementary
  operators to obtain  boundedness from $L^2(\R^n)\times L^2(\R^n)$  to $L^1(\R^n)$. However, even if the linear operator ${\mathscr R}_{\lambda_1}$
  is not well-defined on $L^2(\R^n)$, it is interesting to observe that the average of such operators is well-defined.
  Indeed, a simple computation gives
$$
 \int_{\lambda}^{\mu} {\mathscr R}_{u}(f)(x) du = \int_{\lambda \leq |\xi|\leq \mu}  e^{2\pi i x\cdot \xi }
 {\widehat f}(\xi ) d\xi
 $$
 which is bounded on $L^2(\R^n)$  and one has
 \begin{equation}
 \sup_{\lambda<\mu}
 \left\| \int_{\lambda}^{\mu} {\mathscr R}_{u}  du \right\|_{L^2\to L^2} \leq 1.\label{eq:221}
 \end{equation}
Moreover, in view of the celebrated  result of Fefferman \cite{cfefferman71},
 this operator is unbounded   on $L^p(\R^n)$ if   $p\neq 2$ (as soon as $n\geq 2$).

So we can obtain boundedness from $L^2(\R^n) \times L^2(\R^n)$ to $L^1(\R^n)$ without employing a spherical decomposition but via a
 decomposition along a scale of  ``smoother" operators.
\end{remark}

Following the previous remark, we have the following observation concerning the $L^2\times L^2 \to L^1$ boundedness of
bilinear multipliers.

 \begin{lemma}\label{le3.5} Suppose  that $m_0$
is a bounded function
 with  support in $[-1,1]\times [-1, 1]$ which satisfies
 $$
 \partial_{\lambda_1} \partial_{\lambda_2} m_0(\lambda_1,\lambda_2) \in L^1(\R^2).
 $$
 Let  $m(\xi,\eta)=  m_0(|\xi|,|\eta|)$.
Then $T_{m}$ is bounded from  $L^{2}(\R^n)\times L^{2}(\R^n)$ to $L^1(\R^n)$;
Moreover,
\begin{equation}
\label{eqq} \|T_{m}\|_{L^{2}\times L^{2}\to L^1}\lesssim
  \iint_{\R^2} \left| \partial_{\lambda_1} \partial_{\lambda_2} m_0(\lambda_1,\lambda_2)\right|  d\lambda_1 \, d\lambda_2 .
\end{equation}
\end{lemma}

\medskip

\noindent
{\bf Proof.} We employ a proof via a  decomposition of the symbol as an average of bilinear restriction
operators. First, by modulation and dilation we may assume that $m$ is supported on $[\frac{1}{2},1]\times
[\frac{1}{2},1]$. So via an integration by parts we have
\begin{align*}
T_m(f,g) & = \iint_{[0,1] \times [0,1]} m_0(\lambda_1,\lambda_2) \, {\mathscr R}_{\lambda_1, \lambda_2}(f,g) \, d\lambda_1 \, d\lambda_2 \\
 & = \iint_{[0,1] \times [0,1]}  \partial_{\lambda_1} \partial_{\lambda_2} m_0(\lambda_1,\lambda_2)
  \left(\int_{0}^{\lambda_1} \int_{0}^{\lambda_2} {\mathscr R}_{a, b}(f,g) da\, db \right) d\lambda_1 \, d\lambda_2\\
 & = \iint_{[0,1] \times [0,1]}  \partial_{\lambda_1} \partial_{\lambda_2} m_0(\lambda_1,\lambda_2)
 \left(\int_{0}^{\lambda_1} {\mathscr R}_{a} (f)\, da
 \int_{0}^{\lambda_2} {\mathscr R}_{ b}(g) \,  db \right)  \, d\lambda_1 \, d\lambda_2.
 \end{align*}
Using (\ref{eq:221}) and the H\"older inequality, we deduce
\begin{align*}
\|T_m\|_{L^2 \times L^2\to L^1} & \leq \iint_{[0,1] \times [0,1]} \left|\partial_{\lambda_1}
\partial_{\lambda_2} m_0(\lambda_1,\lambda_2)\right| d\lambda_1 \, d\lambda_2,
  \end{align*}
which concludes the proof. \hfill{} $\Box$

\medskip

Still concerning the boundedness from $L^2(\R^n) \times L^2(\R^n)$ to $L^1(\R^n)$, we have the following result:

\begin{proposition} \label{prop:221}
Let $m_0$ be an even function supported in $[-1,1]^2$ which satisfies the regularity condition:
$$
\sup_{u\in[-1,1]} \| m_0(|u|, |\cdot| )\|_{W^{1+\alpha,1}(\R)} <\infty
$$
for some $\alpha >0$. Then the bilinear operator $T_m$ associated with the symbol $m(\xi,\eta)=m_0(|\xi|,|\eta|)$ is bounded from
$L^2(\R^n)\times L^2(\R^n)$ to $L^1(\R^n)$.
\end{proposition}

\medskip

\noindent
{\bf Proof.} We begin by expressing the operator $T_m$ as follows:
\begin{align*}
 T_m(f,g)(x)& := \int_{\R^{2n}} e^{2 \pi ix \cdot(\xi+\eta)} \widehat{f}(\xi) \widehat{g}(\eta) m_0(|\xi|,|\eta|) d\xi d\eta \\
                  & =  \iint_{[-1,1]^2} m_0(|u|,|v|) {\mathscr R}_{|u|}(f)(x) {\mathscr R}_{|v|}(f)(x).
\end{align*}
The idea is to  express the  function $m_0(|u|,|v|)$ as a tensorial product,  so that   a product of $L^2$-bounded linear operators appears.
So we fix $u\in[-1,1]$ and we examine the function $v \rightarrow m_0(|u|,|v|)$ which is   supported in
$[-1,1]$ and vanishes at the endpoints $\pm 1$. We can expand this
function  in Fourier series
(by considering a  periodic extension on $\R$ of period $2$) and thus we may write for $u,v\in[-1,1]$
$$
m_0(|u|,|v|) = \sum_{k\in{\mathbb Z}} \gamma_k(u) e^{i\pi k v}
$$
with Fourier coefficients
$$
 \gamma_k(u):=  \frac{1}{2} \int_{-1}^1 e^{-i\pi k v} m_0(|u|,|v|) \, dv.
 $$
These coefficients also satisfy the bound for $\alpha\in(0,1)$
\begin{align}
 |\gamma_k(u)| & \lesssim  (1+|k|)^{-1-\alpha} \| m_0(|u|,|\cdot|)\|_{W^{1+\alpha,1}([-1,1])} \label{eq:coeff}
\end{align}
for some $\alpha>0$.
Since we have
\begin{align*}
 T_m(f,g)(x) & = \frac{1}{2}  \sum_{k\in {\mathbb Z}} \iint_{[-1,1]^2}  \gamma_k(u)  {\mathscr R}_{|u|}(f)(x)
 e^{i\pi k v} {\mathscr R}_{|v|}(f)(x)  dudv \\
 & = \frac{1}{2}  \sum_{k\in {\mathbb Z}} \left(\int_{[-1,1]}  \gamma_k(u)  {\mathscr R}_{|u|}(f)(x) du \right)
 \left( \int_{[-1,1]} e^{i\pi k v} {\mathscr R}_{|v|}(f)(x) dv \right),
\end{align*}
we conclude by H\"older inequality and (\ref{eq:coeff}),
\begin{align*}
&\hspace{-1cm} \| T_m(f,g)\|_{1} \\
\lesssim   & \sum_{k\in {\mathbb Z}} (1+|k|)^{-1-\alpha} \left\|\int_{[-1,1]}  (1+|k|)^{1+\alpha} \gamma_k(u)
{\mathscr R}_{|u|}(f) du \right\|_{2} \left\| \int_{[-1,1]} e^{i\pi k v} {\mathscr R}_{|v|}(f)(x) dv \right\|_{2} \\
\lesssim &  \|f\|_{2} \|g\|_{2}.
\end{align*}
Here we used that $\int_{[-1,1]} (1+|k|)^{1+\alpha} \gamma_k(u)  {\mathscr R}_{|u|} du$ is the linear Fourier
multiplier operator associated with the symbol $ (1+|k|)^{1+\alpha} \gamma_k(|\xi|)$ which is uniformly (in $k$)
 bounded on $L^2$, since the symbol is bounded  with respect to $k$ in view of (\ref{eq:coeff}).
\hfill{}$\Box$

\bigskip

  \subsection{An extension  of Theorem~\ref{th3.4}} Using Lemma~\ref{le2.6}, we
  may argue as in  the proof of Theorem~\ref{th3.4} to obtain
    the following result.  The   proof is similar and for brevity is omitted.

  \begin{theorem}\label{th3.5} Let $n\geq 2$ and   $1\leq q_1, q_2 < {2n/(n+1)}$ and
 let $q_1\leq p_1\leq \infty, q_2\leq p_2\leq \infty $ with
 $1/p=1/p_1+1/p_2 $  and $0<p\leq \infty$.  Also assume that
 $$
 {1\over q_1} + {1\over q_2} -{1\over p} \leq \alpha(q_1, q_2),
 $$
 where $\alpha(q_1, q_2)$ is defined in \eqref{e3.10}.
Suppose  that
  $m_0$ is a  bounded function supported in $ [-1,1]\times [-1, 1]$
 such that $m_0 \in W^{\beta,\, 1}(\R^2)$ for some $\beta>n\alpha(q_1, q_2).$
Let  $m(\xi,\eta):=  m_0(|\xi|,|\eta|)$.
Then $T_{m}$ is bounded from  $L^{p_1}(\R^n)\times L^{p_2}(\R^n)$ to $L^p(\R^n)$. In addition,
\begin{eqnarray*}
  \|T_{m}\|_{L^{p_1}\times L^{p_2}\to L^p} \leq
C\|m_0\|_{W^{\beta, 1}({\R^2})}.
\end{eqnarray*}
\end{theorem}

\bigskip

  \section{ Bounds for bilinear Bochner-Riesz means}
\setcounter{equation}{0}

Consider  the bilinear
 Bochner-Riesz means  of order $\delta$ on $\R^n\times \R^n$,
 given by
 %the symbol
%$$
%m(|\xi|, |\eta|)= \Big(1-  {|\xi|^2+|\eta|^2\over R^2}\Big)^{\delta}_+, \ \ \ \  R>0.
%$$
%That is,
\begin{eqnarray}\label{e4.1} \hspace{1cm}
S_R^{\delta} (f, g)(x)  =
 \iint_{|\xi|^2+|\eta|^2\leq R^2} e^{2\pi ix\cdot (\xi+\eta)}
  \bigg(1-  {|\xi|^2+|\eta|^2\over R^2}\bigg)^{\delta} {\widehat f}(\xi) {\widehat g}(\eta)d\xi d\eta.
\end{eqnarray}
In this section, we investigate the range of ${\delta}$ for which the bilinear
 Bochner-Riesz means $S_R^{\delta}  $ are
bounded from  $L^{p_1}(\R^n)\times L^{p_2}(\R^n)$ to $L^p(\R^n)$.
This boundedness holds independently of the parameter $R>0$, so we   take $R=1$ in our
work and for simplicity we write $S^{\delta}$ instead of $S^{\delta}_1.$

We first describe the results in the one-dimensional setting, there Bochner-Riesz multipliers are closely
related to the problem of the disc multiplier.

\begin{theorem}\label{BR-1dim}
 Let $n=1$. The Bochner-Riesz operator is bounded from $L^{p_1}(\R^n)\times L^{p_2}(\R^n)$ into $L^p(\R^n)$
for $1/p=1/p_1 +1/p_2$ in the following situations:

\begin{itemize}
\item[(i)]
  Strict local $L^2$-case:  $2<p_1,p_2,p'<\infty$ and $\delta\geq 0$.

  \item[(ii)]
  Endpoint cases:  $\{p_1,p_2,p'\}=\{2,2,\infty\}$ and $\delta>0$.

\item[(iii)]   Banach triangle case:  $1\leq p_1,p_2,p' \, \le\,  \infty$ and $\delta>0$.
\end{itemize}
\end{theorem}

\noindent
{\bf Proof.}
The first case is a consequence of the positive result for the disc multiplier problem in \cite{GL}. The endpoint can be obtained
 using a discrete spherical decomposition with \cite[Proposition 6.1]{BG2}. The Banach situation follows from similar arguments
 with \cite[Proposition 6.2]{BG2}. Indeed, let us check the point $L^\infty \times L^\infty\to L^\infty$. Consider a smooth
 decomposition of the symbol $(1-|\xi|^2-|\eta|^2)^{\delta}_+=\sum_{\ell\ge 0} 2^{-\delta \ell} m_\ell(\xi,\eta)$ where $m_\ell$
 is supported in a circular neighborhood of the unit circle of approximate distance $2^{-\ell}$ from the circle.
From \cite[Proposition 6.2]{BG2}, we know that
\begin{equation}
 \| T_{m_\ell}(f,g) \|_{L^p \times L^q \to L^\infty} \lesssim 2^{-\ell(\frac{3}{4q}+\frac{1}{2p}) }  \label{eq:dec}
\end{equation}
for any $2\le p,q<\infty$.
\noindent However using integration by parts, it is easy to check that the bilinear kernel $K_{\ell}(x-y,x-z)$ of $T_{m_\ell}$ satisfies
$$ \left|K_\ell (u,v)\right| \le  \frac{C_N\, 2^{-\ell}}{(1+2^{-l}|(u,v)|)^{N}}
$$
for every $N>0$. In this way, (\ref{eq:dec}) can be improved in some off-diagonal estimates as follows: fix $x_0\in\R$ and
define $I:=[x_0-1,x_0+1]$,
\begin{align*}
 \lefteqn{\left|T_{m_\ell}(f,g)(x) \right|} & & \\
 & &   \lesssim 2^{-\ell(\frac{3}{4q}+\frac{1}{2p})} \left(\|f\|_{L^p(2^{M}I)} \|g\|_{L^q(2^{M}I)} +
  \sum_{k_1,k_2\geq M}  2^{-\max\{k_1,k_2\}} \, 2^{\frac{k_1}{p'}+\frac{k_2}{q'}} \|f\|_{L^p(2^{k_1}I)} \|g\|_{L^q(2^{k_2}I)}\right).
\end{align*}
We also conclude that
\begin{align*}
 \left|T_{m_\ell}(f,g)(x) \right| & \lesssim 2^{-\ell(\frac{3}{4q}+\frac{1}{2p})}
  \left(2^{\frac{M}{p}+\frac{M}{q}} + \sum_{k_1,k_2\geq 0}  2^{-N \max\{k_1,k_2\}+(N-1)\ell} \, 2^{k_1+k_2} \right)  \|f\|_{\infty } \|g\|_{\infty} \\
 & \lesssim 2^{-\ell(\frac{3}{4q}+\frac{1}{2p})} \left(2^{\frac{M}{p}+\frac{M}{q}} + 2^{-M(N-2)+(N-1)\ell} \right)  \|f\|_{\infty } \|g\|_{\infty}.
\end{align*}
If we choose $M,N$ such that $M(N-2)=(N-1)\ell$ then we get
$$ \|T_{m_\ell}\|_{L^\infty \times L^\infty \to L^\infty} \lesssim 2^{-\ell(\frac{3}{4q}+\frac{1}{2p}) }
\left(2^{(\frac{1}{p}+\frac{1}{q})\frac{N-1}{N-2} \ell} \right),$$
which holds for every $p,q\in [2,\infty)$.
By taking $p,q$ sufficiently large, we deduce that
$$
\|T_{m_\ell}\|_{L^\infty \times L^\infty \to L^\infty} \lesssim 2^{\rho \ell}
$$ for every $\rho>0$ as small as possible, which concludes the proof by taking $\rho<\delta$.
\hfill{}$\Box$

 \medskip

 These results are optimal in the strict local $L^2$ case, in the endpoint cases, and on the boundary of the Banach triangle. It still
unknown whether boundedness holds in the  limiting case $\delta=0$ in the interior of the Banach triangle minus the local $L^2$ triangle.

We may therefore focus on the higher-dimensional situation.
 First, we have the following proposition.

\begin{proposition} \label{prop4.1}
Let $1\leq p_1, p_2  \leq \infty$ and ${1/ p}={1/p_1}+{1/p_2 }$ with $0<p\leq \infty.$
 Then we have

\begin{itemize}

\smallskip

\item[(i)] If $\delta>n-1/2$, then
$$
\|S^{\delta} (f, g)\|_{{p}}
\leq  C\|f\|_{{p_1}} \|g\|_{{p_2}}.
$$
\item[(ii)] If  $\delta\leq n({1/p}-1) -{1/2}$, i.e.,
$$
p\leq {2n\over 2n+2 \delta +1},
$$
then $S^{\delta} $ is not bounded from $L^{p_1}(\R^n)\times L^{p_2}(\R^n)$ into $L^p(\R^n).$

\item[(iii)] If $\delta\le n \big|\frac{1}{p}-\frac{1}{2}\big| -\frac12$, then
$S^{\delta} $ is unbounded from $L^{p }(\R^n)\times L^{\infty}(\R^n)$ into $L^{p }(\R^n) $, from
$L^{\infty}(\R^n)\times L^{p }(\R^n)$ into $L^{p }(\R^n) $,
 and also
from  $L^{p}(\R^n)\times L^{p' }(\R^n)$ into $L^{1 }(\R^n) $ for any $1\le p \le \infty$.
\end{itemize}
\end{proposition}

\medskip

 \noindent
{\bf Proof.}\  Note that the kernel of the bilinear Bochner-Riesz means $S^{\delta}$ is
$$
K_\delta (x_1,x_2)=c\, {J_{\delta + n}(2\pi |x|) \over  |x|^{\delta + n}}, \ \ \ \ \ x=(x_1,x_2)
$$
and since $ \alpha > n-1/2 $, we have that this satisfies an estimate of the form:
$$
|K_\delta(x_1,x_2)|\lesssim  {1\over (1+|x|)^{\delta + n+1/2}}
$$
by using properties of Bessel functions.
But for such $\delta$ we have $\delta + n+1/2 >2n$, so the kernel satisfies
$$
|K_\delta(x_1,x_2)|\lesssim  {1\over (1+|x_1|)^{n+\epsilon}} {1\over  (1+|x_2|)^{n+\epsilon}},
$$
for some $\epsilon>0$. It follows that the bilinear operator is bounded by a product of two linear operators,
 each of which has a good integrable kernel. So, (i) follows by H\"older's inequality.

We now prove (ii) by using  an argument as in the proof of Proposition 10.2.3 in \cite{Gra}.
  Let $h\in{\mathscr S}(\R^n)$  be a Schwartz function of ${{\mathbb R}^n}$ satisfying that
 \begin{eqnarray*}
 \widehat{h}(\xi)=\left\{
 \begin{array}{ll}
 1, \qquad  |\xi|\leq 2\\[4pt]
 0, \qquad  |\xi|\geq 4.
 \end{array}
 \right.
\end{eqnarray*}
This gives
$$
S^{\delta}(h,h)(x)=\iint_{|\xi|^2+|\eta|^2\leq 1} (1-|\xi|^2-|\eta|^2)^{\delta}
e^{2\pi ix\cdot (\xi+\eta)} \,d\xi d\eta=c\, \frac{J_{n+\delta}(2\pi |(x,x)|)}{|(x,x)|^{n+\delta}}.
$$
Then $S^{\delta}(h,h)(x)$ is a smooth function that is equal to
\begin{eqnarray*}
% \frac{\sqrt{2}\Gamma(\frac{n}{2}+\delta+1)}{\pi^{\frac{n+1}{2}+\delta}} \
c'\, \frac{\cos(2\pi  \sqrt{2} |x|-
\frac{\pi}{2}(n+\delta+\frac{1}{2}))}{ (   \sqrt{2} |x|)^{n+\delta+\frac{1}{2}}}
+O\left(\frac{1}{|x|^{n+\delta+\frac{3}{2}}}\right)
\end{eqnarray*}
 as $|x|\to \infty$.
Then we have
\begin{eqnarray}\label{e4.2}
|S^{\delta}(h,h)(x)|^p\approx
\frac{1}{|x|^{p(n+\delta+\frac{1}{2})}}
+O\left(\frac{1}{|x|^{p(n+\delta+\frac{3}{2})}}\right)
\end{eqnarray}
  for all $|x|$ satisfying
$$
k+\frac{n+\delta}{4}\leq |x|\leq k+\frac{n+\delta}{4}+\frac{1}{4}
$$
for positive large integers $k$.

Now we observe that the error term in (\ref{e4.2}) is of lower order than the main term  at infinity
and thus it does not affect the behavior of $|S^{\delta}(h,h)(x)|^p$.
 So we conclude that $|S^{\delta}(h,h)(x)|^p$ is not integrable when
 $p(n+\delta+1/2)\leq n$,
  i.e. when $p\leq  {2n}/{(2n+2\delta+1)}$.

To prove the first assertion in (iii), we take the $L^\infty$ function to be $1$, and then $S^\delta(f,1)=B^\delta(f$), where
$f$ is the linear Bochner-Riesz operator. So the conclusion follows from the   linear result.
 The second assertion in (iii) is similar. To prove  the third assertion in (iii), by symmetry we may
 assume that $p\le 2$. It   will suffice  to show that the second dual $(S^\delta)^{*2}$ of $S^\delta$
 is unbounded from $L^p\times L^\infty$ to $L^p$. Let $h$ be the Schwartz function in case (ii). Then
  $(S^\delta)^{*2}(h,1)(x) = c \, |x|^{-n/2-\delta} J_{n/2+\delta}(2\pi |x|/\sqrt{2})$ which is not
  an $L^p$ function if $\delta \le n(1/p-1/2)-1/2$.
\hfill{} $\Box$

\medskip

 \subsection{Bilinear Bochner-Riesz means as bilinear multipliers}
 \subsubsection{Main results}
The aim of this section is to prove the following result.

 \begin{theorem}\label{th4.2} Let $n\geq 2$ and   $1\leq p_1, p_2 < {2n/(n+1)}$ and
 $ {1/p}={1/p_1}+ {1/p_2}$
  and $0<p\leq \infty$. Also let $\alpha(p_1,p_2)$ be as in \eqref{e3.10}.
  If  $\delta>n\alpha(p_1, p_2)-1,$  then  the bilinear Bochner-Riesz means operator
   $S^{\delta}$ is  bounded from $L^{p_1}({\R^n})\times
L^{p_2}(\R^n)$ into $L^p(\R^n)$.
%Moreover, we have
% \begin{eqnarray}\label {e4.3}
% \big\|S_R^{\delta} \big\|_{L^{p_1}\times L^{p_2} \to L^p} =
 %\big\|S^{\delta} \big\|_{L^{p_1}\times L^{p_2} \to L^p}  \leq C.
 %\end{eqnarray}
 \end{theorem}

 \medskip

\noindent{\bf Proof.} The proof   is
a consequence of Theorem~\ref{th3.4} and of  Lemma~\ref{le4.3}  proved below.
\hfill{}$\Box$

\newcommand{\f}{\frac}

\begin{lemma}\label{le4.3}
  Let $1\le q<\infty$ and  $\delta=\sigma+i\tau$. If $0<s<\sigma+{1\over q} $,
   then $(1-|x|^2)_+^\delta\in  W^{s,q}(\mathbb{R}^n)$. Moreover, there
  exist constants $C,c>0$ that depend on $n$, $q$, and $s$ such that
  \begin{eqnarray} \label{e4.4}
\big\| (1-|x|^2)_+^\delta\big\|_{W^{s, q}(\mathbb{R}^n)}\leq C\, e^{c|\tau|}
\end{eqnarray}
as long as $ \sigma \le c'$, where $c'$ is a constant.
\end{lemma}

\medskip

 \noindent {\bf Proof.}
To compute the $ W^{s,q}$ norm of $w(x)=(1-|x|^2)_+^\delta$ on $\R^n$,  we argue as follows (see \cite[Theorem 6.3.2]{BL}):
$$
\|w\|_{W^{s,q}} \approx \| w\|_{L^q} + \|w\|_{\dot W^{s,q}},
$$
where $\dot W^{s,q}$ is the homogeneous Sobolev space defined as
$
\|w\|_{\dot W^{s,q}} = \| \Delta^{s/2} w \|_{L^q}\, .
$
We have that $\Delta^{s/2} w$ is a radial function and we can write
\begin{eqnarray*}
\Delta^{s/2} w (x) & =& c \int_{\R^n} \frac{ J_{n/2+\delta}( 2\pi  |\xi|) }{ |\xi|^{n/2+\delta}} |\xi|^s e^{2\pi i x\cdot \xi} d\xi
\\
& = &  C
\int_0^\infty \frac{ J_{n/2+\delta}( 2\pi  r) }{ r^{n/2+\delta}} r^{s+n-1} \int_{S^{n-1}} e^{2\pi i x\cdot r\theta} d\theta \, dr
\\
& = & C'
\int_0^\infty \frac{ J_{n/2+\delta}(  2\pi r) }{ r^{n/2+\delta}} r^{s+n-1} \frac{J_{(n-2)/2} (2\pi r|x|) }{(r|x|)^{(n-2)/2} }   \, dr\\
& = & C''
\int_0^\infty  \tilde{J}_{n/2+\delta}(   r)   r^{s+n-1}\tilde{J}_{(n-2)/2} (r|x|)    \, dr
  =   I,
\end{eqnarray*}
where we set $\tilde{J}_\nu(t) = J_\nu(t) t^{-\nu}$ for $t>0$. Note that we clearly
have that $|\tilde{J}_\nu(t)| \le C_\nu (1+t)^{-\nu-\f12}$
 for all $\nu$ and $t\ge 0$. Moreover, $\tilde{J}_\nu'(t) =-t\tilde{J}_{\nu+1}(t)$ for all $t>0$.
We consider the following cases:

\medskip

\noindent
{\it Case 1: $|x|\le 1/2$.}

\smallskip

 In this case we introduce a smooth cut-off $\psi(r)$ such that $\psi(r)$ is equal to $1$ for $r\ge 3/2$
and $\psi(r)$ vanishes when $r\le 1$. Then $I$ is equal to the sum
$$
  \int_{1}^\infty \tilde{ J}_{n/2+\delta}( r ) r^{s+n-1}   \tilde{J}_{(n-2)/2} (r|x| )   \psi(r)    \, dr
 $$
 $$
 +
  \int_{0}^{3/2} \tilde{ J}_{n/2+\delta}( r)  r^{s+n-1}   \tilde{J}_{(n-2)/2} (r|x| )  (1- \psi(r) )   \, dr \, .
 $$
The second integral is clearly bounded and hence it lies in $L^q(|x|\le 1/2)$. We focus attention on the
first integral. Using properties of the function $\tilde{J}_\nu$ we write
 $$
  \int_{1}^\infty \tilde{ J}_{n/2+\delta}( r ) r^{s+n-1}   \tilde{J}_{(n-2)/2} (r|x| )   \psi(r)    \, dr
  $$
  $$
  =(-1)^k
  \int_{1}^\infty  \Big(\f{d}{rdr }\Big)^k
  \tilde{ J}_{n/2+\delta-k}( r ) r^{s+n-1}   \tilde{J}_{(n-2)/2} (r|x| )   \psi(r)    \, dr
 $$
 Applying a $k$-fold integration by parts we can write the preceding integral as
 $$
   \int_{1}^{\infty}
  \tilde{ J}_{n/2+\delta-k}( r ) \Big(\f{d}{dr } \f 1r\Big)^k   \Big(
  r^{s+n-1}   \tilde{J}_{(n-2)/2} (r|x| )   \psi(r) \Big)   \,   dr.
 $$
If at least one derivative falls on $\psi(r)$, then the integral is easily shown to be bounded. Thus the
worst term appears when no derivative falls on $\psi$. In this case we have
  $$
  \int_{1}^{\infty} \psi(r)\tilde{J}_{\f{n}{2}+\delta-k} (r )   \sum_{\ell=0}^{k} c_\ell
   r^{s+n-1 -2k+2\ell}
  \tilde{J}_{\f{n-2}{2}+\ell} (r|x|)    |x|^{2\ell}        \, dr\, .
 $$
 We examine the $\ell$th term of the sum when $\ell<k$. In this case we split up the integral in the two cases
 $r\ge |x|^{-1}$ and $1\le r\le |x|^{-1}$. In the   case where $r\ge |x|^{-1}$ the integral   contains a factor of $
 r^{s-\sigma-1-k+\ell}$ and this is absolutely convergent since $s-\sigma<1/q\le 1$ and $k-\ell\ge 1$. The term
 overall produces a factor of the form $|x|^{-s+\sigma+k-\f{n-1}{2}}$ which is in $L^q(|x|\le 1/2)$.   In the case
 where  $1\le r\le |x|^{-1}$ one obtains a factor of $|x|^{-s+\sigma+\f{n-1}{2}+k}$ which is also in
 $L^q(|x|\le 1/2)$.

 It remains to consider the case where $\ell=k$. Here we need to show that the term
   $$
 |x|^{2k}  \int_{1}^{\infty} \psi(r)\tilde{J}_{\f{n}{2}+\delta-k} (r )
   r^{s+n-1 }
  \tilde{J}_{\f{n-2}{2}+k} (r|x|)           \, dr\, .
 $$
 is a convergent integral times a positive power of $|x|$. The part of this integral from $1$ to $|x|^{-1}$ is
 bounded by
  $$
   C|x|^{2k}  \int_{1}^{|x|^{-1}} \psi(r)   \f{1}{   r^{\f{n+1}{2}+\delta-k }    }
   r^{s+n-1 }
      \, dr
 $$
which produces a factor of $|x|^{k-c}$, which lies in $L^q(|x|\le 1/2)$. The part of the integral from
$|x|^{-1}$ to $\infty$ is
  $$
  |x|^{2k} \Big[ \int_{|x|^{-1}}^{\infty} \psi(r)   \f{e^{\pm ir}}{   r^{\f{n+1}{2}+\delta-k }    }
   r^{s+n-1 }
  \f{e^{\pm ir|x|}}{  ( r|x|)^{\f{n-1}{2}+k }    }
        \, dr  +  |x|^{-\f{n+1}{2}-k}
 \int_{|x|^{-1}}^{\infty}
  O\Big(
 \f{ r^{s+n-1 }}{   r^{n+2+\delta }   }  \Big)      \, dr\Big]
 $$
 using the asymptotic behavior of the Bessel functions. The second integral converges absolutely while the
 first   integral  contains the phase $ir(\pm 1 \pm |x|)$ which is never vanishing and so it can
be integrated by parts to show that it converges, since $s-\sigma<1/q\le 1$.  At the end one obtains a
factor of $|x|^{k+c}$ which is in $L^q(|x|\le 1/2)$ if $k$ is large.

\medskip

\noindent
{\it Case 2: $|x|\ge 2$. }

\smallskip

In this case we will use again the smooth cut-off $\psi(r)$ which is equal to $1$ for $r\ge 3/2$
and $\psi(r)$ vanishes when $r\le 1$. Then $I$ is equal to the sum
\begin{align}\begin{split}\label{4.100}
 \frac{1}{ |x|^{n+s} } \int_{1}^\infty \tilde{ J}_{n/2+\delta}( r/|x|) r^{s+n-1}   \tilde{J}_{(n-2)/2} (r )   \psi(r)    \, dr \\
 +\,
 \frac{1}{ |x|^{n+s} } \int_{0}^{3/2} \tilde{ J}_{n/2+\delta}( r/|x|)  r^{s+n-1}
   \tilde{J}_{(n-2)/2} (r )  (1- \psi(r) )   \, dr \, .
\end{split}\end{align}
The second integral in \eqref{4.100} is clearly bounded and
since $|x|^{-n-s}$ lies in $L^q(|x|\ge 2)$
the second term in \eqref{4.100} lies in $L^q(|x|\ge 2)$.

 We write the first term in \eqref{4.100} as
 $$
(-1)^k  \f{1}{|x|^{n+s}}    \int_{1}^{\infty}  \tilde{J}_{n/2+\delta}(   r/|x|) \psi(r )  r^{s+n -1 }
\Big(\f{d}{rdr }\Big)^k  \tilde{J}_{\f{n-2}{2}-k}  (r)         \, dr
  $$
for any $k>0$  and by a $k$-fold integration by parts this is equal to
  $$
 \f{1}{|x|^{n+s}}   \int_{1}^{\infty} \Big(\f{d }{dr } \f{1}{r}\Big)^k \Big(
  \psi(r)  \tilde{J}_{n/2+\delta}(   r/|x|)    r^{s+n-1 }  \Big)
  \tilde{J}_{\f{n-2}{2}-k} (r)         \, dr \, .
  $$
 The worst term appears when no derivative falls on $ \psi(r)$.
 In this case we obtain a term of the form
 $$
  \f{1}{|x|^{n+s}}   \int_{1}^{\infty} \psi(r)
  \sum_{\ell=0}^{k} c_\ell \tilde{J}_{\f{n}{2}+\delta+\ell} (r/|x|) \f{1}{|x|^{2\ell}}
   r^{s+n-1 -2k+2\ell}
  \tilde{J}_{\f{n-2}{2}-k} (r)         \, dr
 $$
When $\ell<k$, the $\ell$th term is estimated by
 $$
C |x|^{-n-s-2\ell } \int_{1}^{\infty}r^{s+n-1-2k+2\ell  }
 \f{1}{(1+r/|x|)^{\f{n+1}{2}+\sigma+\ell } } \f{1}{r^{\f{n-1}{2}-k} }dr\, .
 $$
Considering the cases
$r\le |x|$ and $r\ge |x|$ separately, in each case we obtain a convergent integral times a factor of $|x|^{c-k}$,
and since $k$ is arbitrarily large, we deduce that this lies in   $L^q$
  in the range $|x|\ge 2$.  (For the convergence of the integral in the case $r\ge |x|$ we use that
  $s-\sigma<1/q\le 1$.) For the term $\ell=k$ we need the oscillation of the Bessel function
  to show that the integral
\begin{equation}\label{4.101}
  \f{1}{|x|^{n+s+2k}}   \int_{1}^{\infty} \psi(r)  \tilde{J}_{\f{n}{2}+\delta+k} (r/|x|)
   r^{s+n-1 }
  \tilde{J}_{\f{n-2}{2}-k} (r)         \, dr
\end{equation}
 converges.
We split  \eqref{4.101}  as the sum of the term
\begin{equation*}
  \f{1}{|x|^{n+s+2k}}   \int_{1}^{|x|} \psi(r)  \tilde{J}_{\f{n}{2}+\delta+k} (r/|x|)
   r^{s+n-1 }
  \tilde{J}_{\f{n-2}{2}-k} (r)         \, dr,
\end{equation*}
 which is bounded by
  $$
 C \f{1}{|x|^{n+s+2k}}  \int_{1}^{|x|} \psi(r)
   r^{s+n-1 }
 \f{1}{   r^{\f{n-1}{2}-k }    }        \, dr \le C'\, |x|^{-k+c},
 $$
plus the term
$$
  \f{1}{|x|^{n+s+2k}}
 \Big[  \int_{|x|}^{\infty} \psi(r)  \f{e^{\pm i |x|^{-1}r}} {(r/|x|)^{\f{n+1}{2}+\delta+k} }
   r^{s+n-1 }
 \f{e^{\pm ir}}{   r^{\f{n-1}{2}-k }    }     dr  + |x|^{\f{n+3}{2}+\delta+k}
  \int_{|x|}^{\infty}
  O\Big(
 \f{ r^{s+n-1 }}{   r^{n+2+\delta }    }  \Big)      \, dr\Big]\, .
 $$
  Notice that the phase  $ir(\pm 1\pm |x|^{-1})$
 never vanishes, and since $s-\sigma<1/q\le 1$,   an integration by parts yields an absolutely convergent
 integral times a factor of $|x|^{-k+c}$. If $k$ is large, these terms
  have rapid decay at infinity and thus they lie in $L^q(|x|\ge 2)$.

 \bigskip
 Notice that all integrations by parts have produced constants that grow at most like a multiple of
 $1+|s-\sigma+i\tau|^k$ so far.

\medskip

\noindent
{\it Case 3: $1/2\le |x|\le 2$.}

\smallskip

 The part of integral $I$ over the region
  $r\le 2$ is easily shown to be in $L^\infty$ and thus in $L^q$ of the annulus
  $1/2\le |x|\le 2$. It suffices to consider the part of the integral $I$ over the region $r\ge 2$. Here both
  $r$ and $r|x|$ are greater than $1$ and we
  use the asymptotics of the Bessel function to write this part as a sum of terms of the form
$$
C_1 \int_{2}^{\infty}   r^{s-\delta-1 } e^{ i c_1  (|x|+1)r} dr
 + C_2\int_{2}^{\infty}   r^{s-\delta-1 } e^{ i c_2  (|x|-1)r} dr +O\bigg(
 \int_{2}^{\infty}   r^{s-\delta-2 }   dr\bigg)
$$
for some constants $C_1,C_2,c_1,c_2$.
Of these terms the middle one contains a phase that may be vanishing while the other terms
are bounded by constants that grow at most linearly in $|\tau|$.
Recall that $\delta= \sigma +i \tau$. Now define an analytic function of $\delta$ by setting
$$
I_x(\delta)=\int_{2}^{\infty}   r^{s-\delta-1 } e^{   ic_2(|x|- 1)r} dr\, .
$$
 Notice that when $s-\sigma =-\eps_1<0$ we have
$$
| I_x(\delta)| \le   C' \,   \eps_1   ^{-1}  .
$$
 Also notice that when $|x|\neq 1$ and   $s-\sigma =1-\eps_2<1$ we have that
$$
| I_x(\delta)| \le C'' \Big(1 + \frac{ | \tau | }{\eps_2}\Big)\,  |\,  |x|-1|^{-1}.
$$
 But $I_x(\delta)$ is an analytic function of $\delta$ and
 Hirschman's version of the 3-lines lemma (Lemma 1.3.8 in \cite{Gra0})
 gives that
$$
| I_x(\delta)| \le \,  C(     \delta ) |\,  |x|-1|^{-1/q+\eps }
$$
when $s-\sigma = 1/q-\eps <1/q$.
Here
$$
C(\sigma+i\tau) \le  \exp \Big\{\!\frac{\sin (\pi \sigma)}{2}\!
\int_{-\infty}^{\infty} \Big[
\frac{\log  (C'\eps_1^{-1})}{\cosh (\pi  t   )- \cos (\pi \sigma)} +
\frac{\log ( C'' (1+\eps_2^{-1}|t+\tau|))}{\cosh (\pi  t  )+\cos (\pi \sigma)}\Big]dt\! \Big\} \le C_1\, e^{C_2 |\tau|} \, ,
$$
where the last estimate is seen by estimating the logarithm  by a linear term.
But the function $ |\,  |x|-1|^{-1/q+\epsilon }
$ lies in $L^q(1/2\le |x|\le 2)$.
So, we have proved that when
$$
s < \frac{1}{q}+\sigma
$$
we have that $w\in W^{s,q}(\R^n)$.
 \hfill{} $\Box$

\medskip

\subsubsection{An extension  of Theorem~\ref{th4.2}}

Using  Theorem~\ref{th3.5}
and Lemma~\ref{le4.3}, we can obtain
    the following result.

 \begin{theorem}\label{th4.5} Let $n\geq 2$ and   $1\leq q_1, q_2 < {2n/(n+1)}$ and
 let $q_1\leq p_1\leq \infty, q_2\leq p_2\leq \infty $ satisfy
 $1/p=1/p_1+1/p_2 $  and $0<p\leq \infty$.  Also assume that
 $$
 {1\over q_1} + {1\over q_2} -{1\over p} \leq \alpha(q_1, q_2),
 $$
 where $\alpha(q_1, q_2)$ is defined in \eqref{e3.10}.
If $\delta>n\alpha(q_1, q_2)-1,$ then the bilinear Bochner-Riesz means $S^{\delta}$
  is bounded from  $L^{p_1}(\R^n)\times L^{p_2}(\R^n)$ to $L^p(\R^n)$. In addition,
  for some constant $C=C_\delta$  we have
\begin{eqnarray*}
  \|S^{\delta}\|_{L^{p_1}\times L^{p_2}\to L^p} \leq
C.
\end{eqnarray*}
\end{theorem}

 \begin{remark} The preceding result is interesting for $(p_1,p_2)$ near $(1,1)$
as for points away from $(1,1)$ we will obtain better results. Let us show this claim
 by the two examples $(p_1,p_2)=(1, \infty)$ and $(p_1,p_2)=(1, \frac{2n}{n+1})$.

 \begin{itemize}
 \item First let us examine the point $p_1=1$ and $p_2=\infty$. The previous theorem yields
 that if $\delta>{n-1\over 2} + {2n\over n+1},$ then the operator $S^{\delta}$ is bounded from
 $L^{1}(\R^n)\times L^{\infty}(\R^n)$ into $L^1(\R^n)$.

Indeed, we take $1\leq q_2<2n/(n+1)$ such that $1/q_2\in (a_n, b_n)$ (see \eqref{e3.10}) and $q_1=1$,
and so $q_2\leq p_2=\infty$ and $q_1\leq p_1=1$. By \eqref{e3.10}, we have
$$\alpha (q_1, q_2)= {1\over 2} +{2\over n+1} +{1\over 2n}.
$$
On the other hand, ${1\over q_1} +{1\over  q_2} -1\leq \alpha (q_1, q_2). $
By Theorem \ref{th4.5}, we have that if
$ \delta>n \alpha (q_1, q_2) -1= {n-1\over 2} + {2n\over n+1}, $
then the operator $S^{\delta}$ is bounded from $L^{1}(\R^n)\times L^{\infty}(\R^n)$ into $L^1(\R^n)$.

However, we will see in the next section, that for some particular points, such as those with
$p_1=1$ and $p_2=\infty$,
we   have a better result ($\delta>\frac{n}{2}$) using   more precisely the structure of the symbol.

\item Let us now focus on the point $p_1=1$ and $p_2=\frac{2n}{n+1}$. By interpolation between $(1,1,{1\over 2})$
and $(1, \infty, 1)$, Theorem \ref{thm4.88} below and Proposition \ref{prop4.1} imply that $S^\delta$ is bounded
on $L^1(\R^n)\times L^{\frac{2n}{n+1}}(\R^n)$ if $\delta> \frac{3n-2}{4}+\frac{1}{4n}$. In this situation, Theorem \ref{th4.2}
proves that $\delta>n\alpha(\frac{2n}{n+1},1)-1=\frac{n-1}{2}+\frac{2n}{n+1}$ is only necessary (which is better).
\end{itemize}
\end{remark}

\medskip

\subsection{Study of  particular points }

We now focus on determining the range of ${\delta}$
for which the bilinear   Bochner-Riesz means $S^{\delta}  $ are
bounded from  $L^{p_1}(\R^n)\times L^{p_2}(\R^n)$ to $L^{p}(\R^n)$,
when $1/p_1+1/p_2=1/p$, for some specific triples of points $(p_1,p_2,p)$.

\subsubsection{The point $(2,2,1)$ and its dual $(2,\infty,2)$}

 We may easily obtain that  the operator $S^{\delta}$ is bounded from
 $L^2(\R^n)\times L^2(\R^n)$ into $L^1(\R^n)$ when   $\delta>1$. Indeed, to see this
we apply Lemma \ref{le3.5}, so
$$
 \big\|S^{\delta} \big\|_{L^{2}\times L^{2} \to L^{1}}  \lesssim    \iint_{\R^2} \left|\partial_{\lambda_1}
  \partial_{\lambda_2} m(\lambda_1,\lambda_2) \right| d\lambda_1 d\lambda_2,
 $$
where $m(\lambda_1,\lambda_2)=(1-\lambda_1^2-\lambda_2^2)_+^\delta$.
On the disc, we have
$$
 \left|\partial_{\lambda_1} \partial_{\lambda_2} m(\lambda) \right| \lesssim (1-\lambda_1^2-\lambda_2^2)_+^{\delta-2}
 $$
and we can then compute the $L^1$-norm:
\begin{align*}
  \iint_{\R^2} \left|\partial_{\lambda_1} \partial_{\lambda_2} m( \lambda_1,\lambda_2) \right|
 d\lambda_1 d\lambda_2 & \lesssim  \iint_{\R^2} (1-\lambda_1^2-\lambda_2^2)_+^{\delta-2} d\lambda_1 d\lambda_2
    \lesssim 1+\int_{\frac{1}{2}}^1 (1-u)^{\delta-2} du,
  \end{align*}
  which is finite  when $\delta-2>-1$.

The restriction $\delta>1$ is not necessary as shown in our next result:

\begin{theorem} \label{thm:221} Let  $n\geq 2$ and $\delta>0$. Then the operator $S^{\delta}$ is bounded from
 $L^2(\R^n)\times L^2(\R^n)$ to $L^1(\R^n)$. That is, for some constant $C=C_\delta$  we have
 \begin{eqnarray}\label {e4.5bis}
 \big\|S^{\delta} \big\|_{L^{2}\times L^{2} \to L^{1}}    \leq C.
 \end{eqnarray}
 Moreover, this result fails when $\delta=0$.
\end{theorem}

\noindent
{\bf Proof.} As an application of Proposition \ref{prop:221}, we have just to check that for $\delta>0$, there exists $\alpha>0$ with
\begin{equation}\label{4.new100}
\sup_{u\in [-1,1]} \|(1-u^2-\cdot^2)_+^{\delta}\|_{W^{1+\alpha, 1}({\mathbb R})} <\infty.
\end{equation}
Denote by  ${\dot W}^{s, p}({\mathbb R}^d)$ the homogeneous Sobolev space on $\R^d$ with
norm
$$
\|h\|_{{\dot W}^{s, p}({\mathbb R}^d)} = \| (-\Delta)^{s/2} h\|_{L^p(\R^d)}.
$$
 For any $r>0$, we have
$$
\|h(r \cdot)\|_{{\dot W}^{s, p}({\mathbb R}^d)}= r^{-{d\over p} +s} \|h(\cdot)\|_{{\dot W}^{s, p}({\mathbb R}^d)}.
$$
Now for a fixed $\delta>0$ we pick  $0<\alpha<\delta$.   By Lemma \ref{le4.3} we have,
$$
\| (1-|v|^2)_+^{\delta}\|_{{\dot W}^{1+\alpha, 1}({\mathbb R})}\leq C <\infty.
$$
Let   $f(v)=(1-|v|^2)_+^{\delta}$ and $r_u={1\over \sqrt{1-u^2}}$. Combining the preceding facts we obtain
\begin{eqnarray}
\|(1-u^2-\cdot^2)_+^{\delta}\|_{{\dot W}^{1+\alpha, 1}({\mathbb R})}
&=&(1-u^2)^{\delta} \| f( r_u\cdot)\|_{{\dot W}^{1+\alpha, 1}({\mathbb R})}\nonumber\\
&=&(1-u^2)^{\delta-\alpha}  \| f(\cdot)\|_{{\dot W}^{1+\alpha, 1}({\mathbb R})}\nonumber\\
&\lesssim &    (1-u^2)^{\delta-\alpha}\nonumber\\
&\lesssim & 1, \nonumber
\end{eqnarray}
since  $\delta-\alpha>0$  and  then  $(1-u^2)^{\delta-\alpha}\lesssim 1$  for $u\in [-1,1]$.

Certainly, $\|(1-u^2-\cdot^2)_+^{\delta}\|_{L^1({\mathbb R})} <\infty$. Hence (\ref{4.new100}) is proved
and then we conclude the proof by invoking Proposition \ref{prop:221}.

 We now turn to the sharpness of the requirement that $\delta  $ be positive. Let $\mathbb B'$ be the unit ball in
$\R^{2n}$.  If the ball multiplier
$$
T_{\chi_{\mathbb B'}}(f,g) (x)=  \iint_{\R^{2n}} \widehat{f}(\xi) \widehat{g}(\eta) \chi_{\mathbb B'}(\xi,\eta)
e^{2\pi i x\cdot (\xi+\eta)} d\xi d\eta
$$
were bounded from $L^2(\R^n) \times L^2(\R^n)\to L^1(\R^n)$ with norm $C_0$, then by a simple
translation and dilation the multipliers
$T_{\chi_{\mathbb B'_{v,w,\rho}}}$ would also be bounded on the same spaces with norm  $C_0$, where
$$
\mathbb B'_{w,v,\rho} =\{ (\xi,\eta)\in \R^n\times \R^n:\,\, |\xi-\rho w|^2 + |\eta-\rho v|^2 \le 2\rho^2\}.
$$
uniformly for all $\rho>0$ and all unit vectors $v$, $w$ in $\R^n$. Letting $\rho\to \infty$ we would obtain that the
operators
$$
T_{\chi_{P_{w,v}'}}(f,g) (x)=  \iint\limits_{\xi\cdot w+\eta\cdot v\ge 0}\widehat{f}(\xi) \widehat{g}(\eta)
e^{2\pi i x\cdot (\xi+\eta)} d\xi d\eta
$$
are bounded from $L^2(\R^n) \times L^2(\R^n)\to L^1(\R^n)$ with norm $C_0$ uniformly in $v,w\in \mathbb
S^{n-1}$; here $P_{v,w}' =\{(\xi,\eta) \in \R^n\times \R^n:\,\,  \xi\cdot w+\eta\cdot v\ge 0\}$
is a half-space in $\R^{2n}$.
Let $P_{v} =\{\xi \in \R^n:\,\,  \xi\cdot v\ge 0\}$ be a half-space in $\R^n$ determined by $v$.  A simple calculation shows that
$$
T_{\chi_{P_{v,v}'}}(f,g) = ( \widehat{fg} \,\, \chi_{ P_v}     )\spcheck,
$$
and this operator is   unbounded from $L^2\times L^2 \to L^1$  by taking $f=g=\chi_U$, where $U$ is the unit
cube in $\R^n$ and $v=(1,0,\dots , 0)$. This produces a contradiction.
\hfill{}$\Box$

\medskip

 We note that a modification of the preceding counterexample also proves that $S^0$ does not map
  $L^2(\R^n)\times L^\infty(\R^n)$  to $L^2(\R^n)$. Indeed, we take $v $   in $\mathbb S^{n-1}$ and define balls
 $$
\mathbb B''_{ v,\rho} =\{ (\xi,\eta)\in \R^n\times \R^n:\,\, |\xi |^2 + |\eta-\rho v|^2 \le \rho^2\},
$$
 which converge to     $\{ (\xi,\eta)\in \R^n\times \R^n:\,\, \xi\cdot v \ge 0\}$ when $\rho\to \infty$.
 Then  one obtains the operator
 $f  ( \widehat{g} \,\, \chi_{ P_v}     )\spcheck$ in the limit which is unbounded from $L^2(\R^n)\times L^\infty(\R^n)$  to $L^2(\R^n)$.

In the positive direction we show that for $\delta>\frac{n-1}{2}$ boundedness holds in this case.
As of this writing we are uncertain as to whether boundedness holds for the intermediate $\delta$.

 \begin{theorem}\label{th4.7} If $\delta>\frac{n-1}{2}$, then the operator $S^{\delta}$ is bounded from
 $L^2(\R^n)\times L^\infty(\R^n)$  to $L^2(\R^n)$. Moreover, for some constant $C=C_\delta$  we have
 \begin{eqnarray}\label {e4.5}
 \big\|S^{\delta} \big\|_{L^{2}\times L^{\infty} \to L^{2}}    \leq C.
 \end{eqnarray}
\end{theorem}

\medskip

\noindent
{\bf Proof.} Recall that $S^{\delta}$ is a bilinear multiplier with the symbol $
 m_0(|\xi|, |\eta|)= \big(1-  {|\xi|^2-|\eta|^2 }\big)^{\delta}_+.
 $
We now perform a ``spherical decomposition''.
To this end, we choose a smooth function $\chi$ supported on $[1/2,2]$,  which is equal to $1$ on $[3/4, 5/4]$, and which  satisfies
$$
 \sum_{j\geq 0} \chi(2^j t) =1
 $$
for every $t\in(0,1]$. For $j\ge 0$ we introduce the functions
$$
m^j_0(s,t) = (1-s^2-t^2)^{\delta}_+ \, \chi\left(2^j(1-  s^2-t^2) \right).
$$
These symbols give us a spherical decomposition of our initial symbol such that
\begin{equation} \label{eq:dif}
\left|\partial_s^\alpha \partial_t^\beta m_0^j(s,t) \right| \lesssim 2^{-\delta j} 2^{(\alpha+\beta)j}.
\end{equation}
For each such function $m_0^j$, we have the bilinear symbol $m^j(\xi,\eta):=m_0^j(|\xi|,|\eta|)$ on $\R^{2n}$. From (\ref{eq:dif}),
the bilinear operator $T_{m^j}$ has a bilinear kernel $K_{m^j}(x-y,x-z)$ satisfying
$$
\left|K_{m^j}(x-y,x-z)\right| \lesssim 2^{-\delta j} 2^{-j} \left(1+2^{-j}|x-y|+2^{-j}|x-z|\right)^{-M}
$$
for every large enough integer $M>0$.

We may also apply Lemma \ref{le2.6-bis} and from part (i) of Lemma~\ref{le2.2}, we get there exists a constant $C>0$ such that
 \begin{align*}
 \|T_{m^{j}}\|_{L^{2}\times L^{\infty}\to L^2}
  &\leq C   2^{j  (n/2+\epsilon) }
  \| T_{m^{j}}  \|_{L^{2}\times L^{2}\to L^2} + 2^{-N j} \\
   &\leq C   2^{j(n/2+\epsilon) } \sup_{\xi\in\R^n} \Big(\int_{\R^n} |m_0^{j}(|\xi-\eta|, |\eta|)|^2 d\eta\Big)^{1/2} + 2^{-Nj},
 \end{align*}
where $\epsilon>0$ (resp. $N>0$) can be chosen as small (resp. large) as we want.
Moreover, $m_0^{j}$ is supported in the set
$$
  \{(s,t):\ 1-2\cdot 2^{-j } \leq |(s,t)|^2 \leq 1-\tfrac12 \, 2^{-j}\}.
 $$
So from (\ref{eq:dif}), we have
$$
\left(\int_{\R^n} |m_0^{j}(|\xi-\eta|, |\eta|)|^2 d\eta\right)^{1/2} \lesssim 2^{-\delta j} 2^{-\frac{j}{2}},
$$
where we used standard estimates for sub-level sets in dimension $n\geq 2$.
So finally, we conclude that
\begin{align*}
 \|S^{\delta}\|_{L^{2}\times L^{\infty}\to L^2} & \lesssim \sum_{j\geq 0} \|T_{m^j}\|_{L^{2}\times L^{\infty}\to L^2} \\
 & \lesssim \sum_{j\geq 0} \left[  2^{j(n/2+\epsilon) } 2^{-\delta j} 2^{-\frac{j}{2}} + 2^{-N j} \right] \\
 & \lesssim 1 + \sum_{j\geq 0} 2^{-j(\delta -\frac{n-1}{2}-\epsilon)}.
\end{align*}
The proof is then finished since for every $\delta>\frac{n-1}{2}$, we may choose a small enough $\epsilon>0$
such that $\delta> \frac{n-1}{2} + \epsilon$ in order that the previous sum be finite.
\hfill{}$\Box$

\medskip

\subsubsection{The point $(1, \infty, 1)$}

 Related to this point we have the following result which should be contrasted with the unboundedness
 known in this case when $\delta\le \frac{n-1}{2}$ (cf. Proposition \ref{prop4.1} (iii)).

 \begin{theorem}\label{thm4.88}  Suppose  $n\geq 2$.  If $\delta>{n\over 2}$, then the operator $S^{\delta}$ is bounded from
 $L^{1}(\R^n)\times L^{\infty}(\R^n)$ into $L^1(\R^n)$. Moreover, for some constant $C=C_\delta$  we have
 \begin{eqnarray*}
 \big\|S^{\delta} \big\|_{L^{1}\times L^{\infty} \to L^{1}}  \leq C.
 \end{eqnarray*}
\end{theorem}

\noindent
{\bf Proof:}  The proof relies on a mixture of arguments involving in Lemma \ref{le2.6-bis} and the optimal result Theorem \ref{thm:221}.

As previously, express $m_0$ in terms of its   spherical decomposition
$$
 m_0 = \sum_{j\geq 0} m_0^j.
 $$
We have seen (in the proof of Theorem \ref{th4.7}), that we can apply Lemma \ref{le2.6-bis}, which gives
\begin{align*}
 \|S^{\delta}\|_{L^{1}\times L^{\infty}\to L^1} & \lesssim \sum_{j\geq 0} \|T_{m^j}\|_{L^{1}\times L^{\infty}\to L^1} \\
& \lesssim 1+ \sum_{j\geq 0} 2^{j(\frac{n}{2}+\epsilon)}\|T_{m^j}\|_{L^{1}\times L^{2}\to L^1}.
\end{align*}
By Bernstein's inequality (since we only deal with bounded frequencies), it follows  that
\begin{align*}
 \|S^{\delta}\|_{L^{1}\times L^{\infty}\to L^1}  \lesssim 1+ \sum_{j\geq 0} 2^{j(\frac{n}{2}+\epsilon)}\|T_{m^j}\|_{L^{2}\times L^{2}\to L^1}.
\end{align*}
According to Theorem \ref{thm:221} and Proposition \ref{prop:221}, we have
$$
 \|T_{m^j}\|_{L^{2}\times L^{2}\to L^1} \lesssim \sup_{u\in[-1,1]} \| m_0^j(|u|, \cdot )\|_{W^{1+\alpha,1}(\R)}
 $$
for $\alpha >0$ (as small as we wish).
Since
$$ m^j_0(s,t) =  (1-s^2-t^2)^{\delta}_+ \, \chi\left(2^j(1-  s^2-t^2) \right),$$
we have
\begin{align*}
\| m_0^j(|u|, \cdot )\|_{W^{1+\alpha,1}(\R)} & \lesssim 2^{j(1+\alpha-\delta)} 2^{-j} \lesssim 2^{j(\alpha-\delta)}.
\end{align*}
Consequently, we deduce that
\begin{align*}
 \|S^{\delta}\|_{L^{1}\times L^{\infty}\to L^1}  \lesssim 1+ \sum_{j\geq 0} 2^{j(\frac{n}{2}+\epsilon+\alpha-\delta)} <\infty,
\end{align*}
since $\epsilon,\alpha$ can be chosen arbitrarily small and $\delta>\frac{n}{2}$.
\hfill{}$\Box$

 \medskip

\subsection{Interpolation between the different  results}

  Interpolation for $S^\delta$ can be achieved using the bilinear complex method adapted to the setting of
  analytic families  or via the an alternative argument, which is based on bilinear interpolation using the
   real method \cite{4au}. The latter argument is outlined as follows:
We fix $j\geq 0$ and obtain intermediate estimates for each $T_{m^j}$ (depending on $j$) starting from the
 existing estimates for given points.
Since we are still working with ``open conditions'', i.e., a strict inequality of the type $\delta>\delta_0$,
we may obtain intermediate boundedness for $S^\delta$   by interpolating the boundary conditions.

\begin{comment}
 \begin{definition}\label{def4.10} let $T$ be a linear operator bounded from
 $L^{p_1}(\R^n)\times L^{p_2}(\R^n)$ into $L^p(\R^n)$ for exponents $p, p_1, p_2\in [1, \infty].$
 Using real duality, we define its two adjoints $T^{\ast 1}$ and $T^{\ast 2}$ by
 $$
 \langle T(f,g), h\rangle =: \langle T^{\ast 1}(h,g), f\rangle =:\langle T^{\ast 2}(f,h), g\rangle.
 $$
 So $T^{\ast 1}$ is bounded from $L^{p'}(\R^n)\times L^{p_2}(\R^n)$ into $L^{p'_1}(\R^n)$  and $T^{\ast 2}$
 is bounded from $L^{p_1}(\R^n)\times L^{p'}(\R^n)$ into $L^{p'_2}(\R^n).$
 \end{definition}

 \medskip

 Having obtained the above bilinear estimates, we can interpolate between them. A simple computation gives that
  $$
T^{\ast 1}(h,g)(x) =  \int_{\R^{2n}} e^{2\pi ix\cdot (\xi+\eta)} m(-\xi-\eta, \xi)
 {\widehat h}(\xi) {\widehat g}(\eta) d\xi d\eta.
 $$
 So  it corresponds to the bilinear integral associated with the symbol
 $$
 m^{\ast 1}(\xi, \eta) =m(-\xi-\eta, \eta).
 $$
  Similarly $T^{\ast 2}$  corresponds to the bilinear integral associated with the symbol
 $$
 m^{\ast 2}(\xi, \eta)=m(\xi, -\xi-\eta).
 $$
\end{comment}

\begin{center}

\hspace{.5in} \scalebox{1.50}{\includegraphics{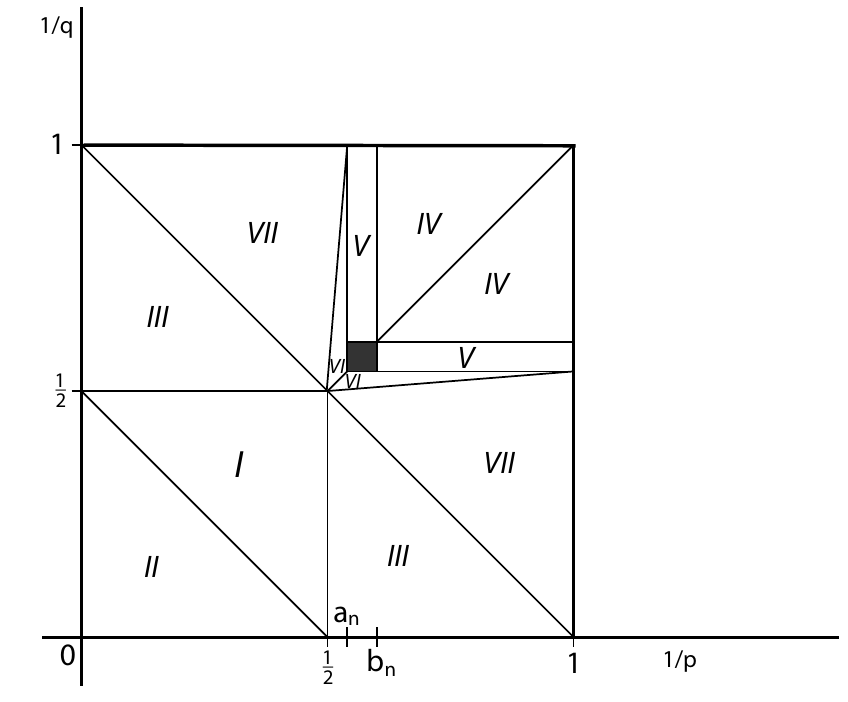} }

\hspace{-.1in}
{\small  Figure : Exponents $(\frac{1}{p},\frac{1}{q})$ for $p,q\geq 1$. Here $a_n={n+1\over 2n}$, $b_n={n+1\over 2n} +{n-1\over n^2+n}$}

\end{center}

\bigskip

\bigskip

We consider a spherical decomposition, as   in Theorem \ref{th4.7} splitting the symbol $m_\delta(\xi,\eta) =(1-|(\xi,\eta)|)_+^\delta$  as
\begin{equation} \label{eq:deco}
m_\delta = \sum_{j\geq 0} 2^{-j \delta}  m^{j,\delta}
\end{equation}
with bi-radial symbols  $m^{j,\delta}(\xi,\eta)= m_0^{j,\delta}(|\xi|,|\eta|)$, where
$$
m_0^{j,\delta}(t,s)= \big[2^j (1-t^2-s^2)\big]^\delta \chi \big(2^j(1-t^2-s^2)\big)  ,
$$
and thus $m^{j,\delta}$ are
supported in the annulus $1-|(\xi,\eta)| \simeq 2^{-j}$ and are regular at the scale $2^{-j}.$

In the proof of boundedness of  $T_{m_\delta}$  at the points
 $$
 (p_1,p_2,p)\in\mathcal S=\{(2,2,1),(1,\infty,1),(\infty, 1,1), (2,\infty,2), (\infty, 2,2) , \dots \}  ,
 $$
  we actually obtained estimates of the form

\begin{equation}\label{bbgghh} \|  T_{m^{j,\delta}} \|_{L^{p_1} \times L^{p_2} \to L^p} \lesssim 2^{j\delta(p_1,p_2,p)}
\end{equation}
for some $\delta(p_1,p_2,p)$ depending only on the  points $p_1,p_2,p$.
In proving \eqref{bbgghh}, we only used the biradial nature of $m^{j,\delta}$, its support properties, and the bounds
$$
| \partial_t^\alpha \partial_s^\beta m_0^{j,\delta} (t,s)| \lesssim  2^{(\alpha+\beta)j}
$$
   that are independent of  $\delta$;   thus estimate \eqref{bbgghh} also holds for any other $m^{j,\delta'}$, i.e.,
   \begin{equation}\label{bbgghh2} \|  T_{m^{j,\delta'}} \|_{L^{p_1} \times L^{p_2} \to L^p} \lesssim 2^{j\delta(p_1,p_2,p)}.
\end{equation}

   We now  fix $j$ and $\delta'$ and apply estimate \eqref{bbgghh2} and
   bilinear real interpolation  (as in \cite{4au})  on $T_{m^{p,\delta'}}$ between
   the points $(p_1^0,p_2^0,p^0)$ and $(p_1^1,p_2^1,p^1)$ to obtain a bound
   $$
   \|  T_{m^{j,\delta'}} \|_{L^{p_1} \times L^{p_2} \to L^p} \lesssim (2^{j})^{(1-\theta)\delta(p_1^0,p_2^0,p^0)+
   \theta \delta(p_1^1,p_2^1,p^1)}.
   $$
   Define $\delta(p_1,p_2,p)=(1-\theta)\delta(p_1^0,p_2^0,p^0)+
   \theta \delta(p_1^1,p_2^1,p^1)$ whenever $(p_1^0,p_2^0,p^0), (p_1^1,p_2^1,p^1)$ are    in $\mathcal S$
   and
   $$
  \Big(\f{1}{p_1},\f{1}{p_2}, \f{1}{p}\Big) =(1-\theta)   \Big(\f{1}{p_1^0},\f{1}{p_2^0}, \f{1}{p^0}\Big)+ \theta
  \Big(\f{1}{p_1^1},\f{1}{p_2^1}, \f{1}{p^1}\Big).
   $$
Then we obtain the bound \eqref{bbgghh2} for  $T_{m^{j,\delta'}}$  and  for any triple of points  $(p_1,p_2,p)$ and any $\delta'$.
 Picking $\delta'=\delta$ and summing over $j$ yields a bound for $S^\delta$ from
 $L^{p_1}\times L^{p_2}\to L^p$ when $\delta>\delta(p_1,p_2,p)$.  The summation over $j$ is straightforward
  when $p\ge 1$. In the case where $p\le 1$ we sum the series  as follows:
 $$
\Big\| \sum_{j\ge 0 } 2^{-j\delta} T_{m^{j,\delta}} \Big\|_{L^{p_1} \times L^{p_2} \to L^p }^p\le
\sum_{j\ge 0 } 2^{-j\delta p} \big\| T_{m^{j,\delta}} \big\|_{L^{p_1} \times L^{p_2} \to L^p }^p
\lesssim \sum_{j\ge 0 } 2^{-j\delta p} (2^j)^{p\delta(p_1,p_2,p)}<\infty\, ,
 $$
 which also converges as long as $\delta>\delta(p_1,p_2,p)$.

\bigskip

Via this method we obtain the following results:

\begin{proposition}[Local-$L^2$ case] Let  $p,q,r'\in[2,\infty)$ with $\frac{1}{p}+\frac{1}{q}+\frac{1}{r'}=1$ (region $I$).
If $\delta> \frac{n-1}{r'}$, then $S^\delta$ is bounded from $L^p(\R^n) \times L^q(\R^n)$ into $L^r(\R^n)$.
\end{proposition}

\begin{proposition}[Banach case]
\begin{itemize}
\item[(a)]   Let $p,q\in[2,\infty)$ and $r'<2$  with $\frac{1}{p}+\frac{1}{q}+\frac{1}{r'}=1$ (region $II$). If $\delta> \frac{n-1}{2} + n(\frac{1}{r'}-\frac{1}{2})$, then
   $S^\delta$ is bounded from $L^p(\R^n) \times L^q(\R^n)$ to $L^r(\R^n)$.

   \smallskip

\item[(b)]  Let $q,r'\in[2,\infty)$ and $p<2$  with $\frac{1}{p}+\frac{1}{q}+\frac{1}{r'}=1$ (region $III$). If $\delta> n(\frac{1}{2}-\frac{1}{q}) -\frac{1}{r'}$, then
  $S^\delta$ is bounded from $L^p(\R^n) \times L^q(\R^n)$ to $L^r(\R^n)$.

   \smallskip
\item[(c)]   Let $p,r'\in[2,\infty)$ and $q<2$ with $\frac{1}{p}+\frac{1}{q}+\frac{1}{r'}=1$. If  $\delta> n(\frac{1}{2}-\frac{1}{p}) -\frac{1}{r'}$, then
   $S^\delta$ is bounded from $L^p(\R^n) \times L^q(\R^n)$ to $L^r(\R^n)$.

  \end{itemize}
\end{proposition}

We now address the  non-Banach case situation which is more complicated: if $q\geq p$, then interpolating
between the point $(1,1,\frac{1}{2})$  and Theorems \ref{thm:221} and \ref{thm4.88}
yields
$$ \delta > \delta_1:=n\Big(\frac{1}{p}-\frac{1}{2}\Big)-\frac{n-1}{2r'}.$$
(Here we recall that $\frac{1}{r'}=1-\frac{1}{r}\leq 0$).
But Theorem \ref{th4.2} (for $a_n \leq \frac{1}{q} \leq \frac{1}{p}\leq 1$)  gives the condition
$$ \delta > \delta_2:= n \alpha(p,q) -1.$$

If $a_n \leq \frac{1}{q} \leq \frac{1}{p} \leq b_n$, then $\alpha(p,q)=\frac{4}{n+1}$ and we check
that $\delta_2\geq \delta_1$, so Theorem \ref{th4.2} does not improve the exponent $\delta_1$ (and the same if $p,q$ are bigger).

If $a_n \leq \frac{1}{q} \leq b_n \leq \frac{1}{p} \leq 1$, then we see that $\delta_2 \leq \delta_1$ if an only if
$$ \frac{1}{r} \geq \frac{3n-1}{n^2-1}+1.$$

If $b_n \leq \frac{1}{q} \leq \frac{1}{p} \leq 1$, then we have $\delta_2 \leq \delta_1$ if and only if
$$ \frac{1}{q} \leq \frac{1}{p} + \frac{1}{nr'}.$$

Collecting this information together, we   deduce the following result.

\begin{proposition}[Non-Banach case, part 1]
Let $p\leq \min\{2,q\}$ and $\frac{1}{r}:=\frac{1}{p}+\frac{1}{q} > 1$. Then $S^\delta$ is bounded
from $L^{p}(\R^n) \times L^q(\R^n)$ into $L^r(\R^n)$
\begin{itemize}

 \item if $b_n \leq \frac{1}{q} \leq \frac{1}{p} \leq 1$ (region $IV$): $\delta>\delta_2$
 for $\frac{1}{q} \leq \frac{1}{p} + \frac{1}{nr'}$, and   $\delta>\delta_1$
  for $\frac{1}{q} > \frac{1}{p} + \frac{1}{nr'}$;

 \item or if $a_n \leq \frac{1}{q} \leq b_n \leq \frac{1}{p} \leq 1$ (region $V$): $\delta>\delta_2$
 for $ \frac{1}{r} \geq \frac{3n-1}{n^2-1}+1$, and $\delta>\delta_1$ for $ \frac{1}{r} < \frac{3n-1}{n^2-1}+1$;

 \item or if $a_n\leq \frac{1}{q} \leq \frac{1}{p} \leq b_n$ (shaded region) and $\delta> \delta_1$.
\end{itemize}
\end{proposition}

We are left with regions $VI$ and $VII$.  For region $VI$ we interpolate between the point $(2,2,1)$
and the line segments $\{(a_n,1/q):\,\, 1/q \in (a_n,1]\}$ and $\{( 1/p,b_n):\,\, 1/p \in (a_n,1]\}$ to obtain the following result:

\begin{proposition}[Non-Banach case, part 2]
For a point $(1/p^0,1/q^0)$ in  the part of region $VI$ above the diagonal
find $\theta\in (0,1)$ and find $1/q\in (a_n,1]$ such that
$1/p^0=(1-\theta)/2+\theta a_n$ and $1/q^0=(1-\theta)/2+\theta /q$. Then $S^\delta$ is bounded from
$L^{p^0}\times L^{q^0}$ to $L^{r^0}$ where $1/r^0=1/p^0+1/q^0$ whenever $\delta>\theta n \alpha(\f{2n}{n+1}, q)-1$.

An analogous result holds for the part of region $VI$ below the diagonal.
 \end{proposition}

 Finally in region $VII$ we apply a similar interpolation between the point $(2,2,1)$ and the line segments
 joining the points $(1,\infty, 1)$ with $(1,a_n, a_n/(a_n+1))$ and
 $( \infty, 1,1)$ with $( a_n, 1,a_n/(a_n+1))$ to obtain $L^{p_1}\times L^{p_2}\to L^p$
 boundedness for $\delta$ bigger than some critical
 value $\delta(p_1,p_2,p)$.

\section{Concluding remarks}

The linear Bochner-Riesz problem has been studied by several authors; we
% cite for instance
refer readers to \cite{Bo, B1, B2, BoG, CS, F1, cfefferman71, F3,
Gra, Lee, Sog,  St2, SW, T2,  TVV} and the references therein for further relevant literature.
We are not sure how to adapt the  techniques in these articles  to the bilinear setting but we hope to investigate whether the bilinear approach to the
restriction and Kakeya conjectures in \cite{TVV} could potentially shed some new light in this problem.

\vskip 1cm

\noindent
{\bf Acknowledgements:}
F. Bernicot  is supported by the ANR under the project AFoMEN no. 2011-JS01-001-01.
L. Grafakos is  supported by  the
National Science Foundation (USA) under grant number 0900946.   L. Song is supported by  NNSF of China (No. 11001276).   L. Yan
is  supported by  NNSF of China (Grant No.  10925106).
L. Yan  would like
 to thank A. Seeger and C. Sogge   for fruitful discussions.

 \vskip 1cm

{\scriptsize { CNRS - Universit\'e de Nantes, Laboratoire Jean Leray  2, rue de la Houssini\`ere
44322 Nantes cedex 3, France

{\it E-mail address}:  frederic.bernicot@univ-nantes.fr}

\medskip

{\scriptsize { Department of Mathematics, University of Misouri, Columbia, MO 65211, USA

{\it E-mail address}:  grafakosl@missouri.edu}

 \medskip

 {\scriptsize {Department of Mathematics, Sun Yat-sen University, Guangzhou, 510275,
 P.R. China

{\it E-mail address}:  songl@mail.sysu.edu.cn

 \medskip

{\scriptsize {
Department of Mathematics, Sun Yat-sen University, Guangzhou, 510275,
 P.R. China

{\it E-mail address}:  mcsylx@mail.sysu.edu.cn
}}}

\end{document}